\documentclass[thmsa,notitlepage,leqno]{article}%
\usepackage{amssymb}
\usepackage{amsfonts}
\usepackage[leqno]{amsmath}
\usepackage{sw20jart}
\usepackage{graphicx}
\usepackage{epstopdf}%
\providecommand{\U}[1]{\protect\rule{.1in}{.1in}}
\begin{document}

\author{Pablo Azcue\thanks{Departamento de Matematicas, Universidad Torcuato Di Tella.
Av. Figueroa Alcorta 7350 (C1428BIJ) Ciudad de Buenos Aires, Argentina.}, Nora
Muler$^{\ast}$ and Zbigniew Palmowski\thanks{Faculty of Pure and Applied
Mathematics, Wroc\l aw University of Science and Technology, Wyb.
Wyspia\'nskiego 27, 50-370 Wroc\l aw, Poland.}}
\title{Optimal dividend payments for a two-dimensional insurance risk process}
\maketitle

\begin{abstract}
We consider a two-dimensional optimal dividend problem in the context of two
branches of an insurance company with compound Poisson surplus processes
dividing claims and premia in some specified proportions. We solve the
stochastic control problem of maximizing expected cumulative discounted
dividend payments (among all admissible dividend strategies) until ruin of at
least one company. We prove that the value function is the smallest viscosity
supersolution of the respective Hamilton-Jacobi-Bellman equation and we
describe the optimal strategy. We analize some numerical examples.

\end{abstract}

\section{Introduction}

In collective risk theory the surplus process $X$ of an insurance company is
modeled as
\begin{equation}
X(t)=u+ct-S(t), \label{mdim}%
\end{equation}
where $u>0$ denotes the initial surplus,
\begin{equation}
S(t)=\sum_{i=1}^{N_{t}}U_{i} \label{zlpoisson}%
\end{equation}
is a compound Poisson. We assume that $U_{i},(i=1,2,...)$ are i.i.d.
distributed claims (with the distribution function $G$).The arrival process is
a homogeneous Poisson process $N_{t}$ with intensity $\lambda$. The premium
income is modeled by a constant premium density $c$ and often assumed net
profit condition $c>\lambda E[U_{1}]$ gives the unrealistic property that
process $X(t)$ converges to infinity. In answer to this objection De Finetti
\cite{Fin} introduced the dividend barrier model for one-dimensional model
(\ref{mdim}), in which all surpluses above a given level are transferred to a
beneficiary. Further, usually the payment of dividends should be made in such
a way as to maximize the expected discounted sum of dividends paid up to ruin
time. In 1969, Gerber \cite{[15]} showed that if the free surplus of an
insurance portfolio is modelled by a compound Poisson risk model, it is always
optimal to pay dividends according to a so-called band strategy, which
collapses to a barrier strategy for exponentially distributed claim amounts.
Later, lots of works in the mathematical finance and actuarial literature
concern the dividend barrier models and the problem of finding an optimal
policy for paying out dividends. Gerber and Shiu \cite{GerberEllias}, Grandits
et al. \cite{GHS} and Jeanblanc and Shiryaev \cite{JeanShir} consider the
optimal dividend problem in a Brownian setting. Irb\"{a}ck \cite{Irback} and
Zhou \cite{Zhou}, Zajic \cite{Zajic}, Avram et al. \cite{APPdiv}, Kyprianou
and Palmowski \cite{KyprPalm}, Loeffen \cite{Loeffen1} study the constant
barrier model for a classical and spectrally negative L\'{e}vy risk process.
Azcue and Muler \cite{azmuazul} follow a viscosity approach to investigate
optimal reinsurance and dividend policies in the Cram\'{e}r-Lundberg model.
The most general criteria currently available for barrier strategies to be
optimal can be found in Loeffen and Renaud \cite{[19]}.

A detailed overviews on this subject from different perspectives are given in
Azcue and Muler \cite{azmu}, Schmidli \cite{[21]}, Albrecher and Thonhauser
\cite{[2]} and Avanzi \cite{[4]}.

All these control problems have been formulated and studied in the
one-dimensional framework. However, in recent years there has been an
increased interest in risk theory in considering the multidimensional surplus
model where $X(t), x, c$ and $S(t)$ are vectors, with possible dependence
between the components. Indeed, the assumption of independence of risks may
easily fail, for example in the case of reinsurance, when incoming claims have
an impact on both insuring companies at the same time. In general, one can
also consider situations where each claim event might induce more than one
type of claim in an umbrella policy (see Sundt \cite{Su}). For some recent
papers considering dependent risks, see Dhaene and Goovaerts \cite{DG, DG2},
Goovaerts and Dhaene \cite{DG3}, M\"{u}ller \cite{Mulla, Mullb}, Denuit et al.
\cite{DGM}, Ambagaspitiya \cite{Am}, Dhaene and Denuit \cite{DD}, Hu and Wu
\cite{HW} and Chan et al. \cite{chan}. Ruin probability expressions for a
two-dimensional risk process were also studied in Avram et al. \cite{[5], [6]}
for simultaneous claim arrivals and proportional claim sizes and recently in
Badila et al. \cite{[11]} and Ivanovs and Boxma \cite{[17]} in a more general framework.

In this paper we analyze the dividend problem for the two-dimensional risk
model in which two branches of a company split the amount they pay out of each
claim in fixed proportions $b_{1}$ and $b_{2}$ ($b_{1}+b_{2}=1$), and receive
premiums at rates $c_{1}$ and $c_{2}$, respectively. Moreover, these two
branches have the same shareholders and the objective is to maximize the total
dividends received by these shareholders. That is, the main goal of this paper
is identification of the value function maximizing the weighted sum of
expected discounted dividend payments until ruin of at least one branch. This
will lead to a fully two-dimensional stochastic control problem and we answer
the question what the analogues of the optimal univariate barrier strategies
are in two dimensions.

In Azcue and Muler \cite{[8]}, the problem of optimally transferring capital
between two portfolios in the presence of transaction costs was considered,
see also Badescu et al. \cite{[10]}. Albrecher et al. \cite{alazmu} considered
the problem of optimizing dividends for two collaborating insurance companies
whose surpluses are modelled with independent compound Poisson processes.

Czarna and Palmowski \cite{CzarPal} studied the same dividend problem as in
this paper but for a very particular dividend strategy of reflecting
two-dimensional risk process from the line. Solving certain partial
differential equations they managed to identify the value function for the
exponential claim sizes. We will show in this paper though that this strategy
is not optimal.

In this paper we prove that the value function is a viscosity solution of the
respective Hamilton-Jacobi-Bellman equation (abbreviated lately by HJB) and it
can be characterized as the smallest viscosity supersolution. Using this
results we manage to describe the optimal strategy for the symmetric case
$c_{2}/c_{1}=b_{2}/b_{1}$. In the remainder cases, we describe that the best
strategy for initial surpluses $(x_{1},x_{2})$ satisfying $x_{1}/b_{1}%
>x_{2}/b_{2}$ is the first branch paying $\left(  x_{1}-\left(  b_{1}%
/b_{2}\right)  x_{2}\right)  $ immediately as dividends. We use a convergent
numerical scheme to find the optimal strategy for initial surpluses
$(x_{1},x_{2})$ satisfying $x_{1}/b_{1}\leq x_{2}/b_{2}$; this numerical
scheme is a particular case of the one presented in Azcue and Muler
\cite{AzMuNum}. We believe that the techniques used in this paper could be
applicable to other risk-sharing mechanisms.

The rest of the paper is organized as follows. In Section \ref{sec:model} we
present the model we deal with and state formally the problem that we want to
solve. In Section \ref{sec:hjb} we prove that the value function is a
viscosity solution of HJB. Using this result we describe the optimal strategy
in some cases in Section \ref{sec: optimalstrategy}, in Section
\ref{NumericalScheme} we describe a numerical scheme to approximate the
optimal value function and the optimal strategy and in Section
\ref{sec:examples}\ we present some numerical examples.
%We conclude
%the paper with the number of examples given in Section \ref{sec:examples}.

\section{Model}

\label{sec:model}

In this paper we consider two branches of a company that receive premiums at
different rates and then split the amount they pay in fixed proportions for
each claim. This model corresponds to proportional reinsurance dependence. The
total amount of claims up to time $t$ is modeled as a compound Poisson
process. We can write this two-dimensional risk model as%

\begin{equation}
\overline{X}(t)=(X_{1}(t),X_{2}(t))=(x_{1}+c_{1}t-%
%TCIMACRO{\dsum \limits_{i=1}^{N_{t}}}%
%BeginExpansion
{\displaystyle\sum\limits_{i=1}^{N_{t}}}
%EndExpansion
b_{1}U_{i},x_{2}+c_{2}t-%
%TCIMACRO{\dsum \limits_{i=1}^{N_{t}}}%
%BeginExpansion
{\displaystyle\sum\limits_{i=1}^{N_{t}}}
%EndExpansion
b_{2}U_{i}). \label{uncontrolledSurplus}%
\end{equation}
Above $x_{1}$\ and $x_{2}$ are the corresponding initial surplus levels and
$c_{1}$ and $c_{2}$ are the corresponding premium rates. The sizes of the
claims $U_{i}\ $are non-negative i.i.d. random variables with absolutely
continuous distribution $G$. The claims arrival process $N_{t}\ $is a Poisson
processes with intensity $\lambda$. Finally, the constants $b_{1}$ and $b_{2}$
are the proportions of the claim that each branch pays, so $b_{1}+b_{2}=1$,
$b_{1}>0$ and $b_{2}>0$. We assume here that the process $N_{t}$\ and the
random variables $U_{i}$ are independent of each other. Without loss of
generality, we can assume that the second branch receives less premium per
amount paid out, that is
\begin{equation}
c_{1}/b_{1}\geq c_{2}/b_{2}. \label{case}%
\end{equation}
We denote by $\{\mathcal{F}_{t}\}_{\{t\geq0\}}$ the right-continuous natural
filtration of $\overline{X}(t)$ satisfying usual conditions. Throughout this
paper all stopping times and martingales are taken with respect of this filtration.

Both branches use part of their surpluses to pay dividends. The dividend
payment strategy $\overline{L}(t)=\left(  L_{1}(t){\small ,}L_{2}(t)\right)
$\ is the total amount of dividends paid by the two branches up to time $t.$
We define the associated controlled process with initial surplus $(x_{1}%
,x_{2})$ as%

\begin{equation}
\overline{X}^{\overline{L}}(t)=\overline{X}(t)-\overline{L}(t).
\label{X_Controlado}%
\end{equation}
The dividend payment strategy $\overline{L}=\left(  L_{1}(t){\small ,}%
L_{2}(t)\right)  _{t\leq\tau}$ is called \textit{admissible} if it is
non-decreasing, c\`{a}gl\`{a}d (left continuous with right limits),
predictable with respect to the filtration generated by the bivariate process
$\overline{X}(t)$ and satisfies $L_{1}(t)\leq X_{1}(t)$, $L_{2}(t)\leq
X_{2}(t)$.
%Let us define $\tau_{i}$ as the arrival time of the $i$-th claim
The first time when the two-dimensional risk process first leaves the positive
quadrant will be our ruin time:
\[
\widehat{\tau}=\inf\left\{  t:X_{1}(t)<0\text{ or }X_{2}(t)<0\right\}  .
\]
That is, the ruin time is the first time at which at least one branch get
ruined. We denote by $\mathbf{R}_{+}=[0,\infty)$ and by $\mathbf{R}_{+}%
^{2}=(0,\infty)^{2}$ the first quadrant. Let $\Pi_{x_{_{1}},x_{2}}$ be the set
of admissible dividend strategies with the initial surplus $\overline
{x}=(x_{1},x_{2})\in\mathbf{R}_{+}^{2}$. Given an admissible dividend strategy
$\overline{L}\in\Pi_{x_{_{1}},x_{2}},$ the expected discounted dividends paid
by the two branches until $\widehat{\tau}$ is%
\[
V_{\overline{L}}(\overline{x})=E_{\overline{x}}\left(  \int_{0}^{\widehat
{\tau}}e^{-qs}dL_{1}(s)+\int_{0}^{\widehat{\tau}}e^{-qs}dL_{2}(s)\right)  ,
\]
where $q>0$ is a constant discount factor.

The main goal of this paper is identification of the optimal value function
defined by:
\begin{equation}
V(\overline{x})=\sup_{\overline{L}\in\Pi_{\overline{x}}}V_{\overline{L}%
}(\overline{x}) \label{V}%
\end{equation}
for $\overline{x}\in\mathbf{R}_{+}^{2}.$

In the next section, we will see that $V$ is well defined. Another crucial
problem that we treated in this paper is the existence of an optimal strategy
defined as follows. Given a family of admissible strategies
\[
\pi=\left\{  \overline{L}_{\overline{x}}\in\Pi_{\overline{x}}\text{ for each
}\overline{x}\in\mathbf{R}_{+}^{2}\right\}
\]
we define the \textit{value function} $V_{\pi}:\mathbf{R}_{+}^{2}%
\rightarrow\mathbf{R}_{+}$ as $V_{\pi}(\overline{x})=V_{\overline
{L}_{\overline{x}}}(\overline{x})$. We say that $\pi^{\ast}$ is the
\textit{optimal strategy} if $V_{\pi^{\ast}}=V$.

\begin{remark}
\label{RemarkComparacionMerger} We can also consider the possibility that the
two branches of the company can merge with a merger cost $m$ (see e.g. Gerber
and Shiu \cite{GerberShiu}). The merged company has initial surplus
$x_{1}+x_{2}-m,$ receives the premium rate $c_{1}+c_{2}$ and pays the whole
claims $U_{i}.$ The new company uses part of the surplus to pay dividends to
the shareholders up to the time in which the joined surplus becomes negative.
In this case, the uncontrolled surplus process is given by%
\[
X_{M}(t)=x_{1}+x_{2}-m+(c_{1}+c_{2})t-%
%TCIMACRO{\dsum \limits_{i=1}^{N_{t}}}%
%BeginExpansion
{\displaystyle\sum\limits_{i=1}^{N_{t}}}
%EndExpansion
U_{i},
\]
and the controlled surplus process is given by $X_{M}^{L}(t)=X_{M}%
(t)-L(t),$where the dividend payment strategy $L(t)$ is non-decreasing,
c\`{a}gl\`{a}d (left continuous with right limits), predictable with respect
to the filtration generated by the process $X_{M}(t)$ and satisfies $L(t)\leq
X_{M}(t)$. We consider the ruin time of the merged company $\tau_{M}^{L}%
=\inf\left\{  t:X_{M}^{L}(t)<0\text{ }\right\}  $ and we define the merger
optimal value function as%
\[
V_{M}(x_{1}+x_{2}-m)=\sup_{L}E_{\overline{x}}\left(  \int_{0}^{\tau_{M}^{L}%
}e^{-qs}dL(s)\right)  .
\]
for $x_{1}+x_{2}\geq m.$ The function $V_{M}$ is the optimal value function of
the one-dimensional De Finetti's problem corresponding to the compound Poisson
process $X_{M}$.

Let us assume that $m=0,$ then given any $\overline{x}\in R_{+}^{2}$ and any
admissible strategy $\overline{L}(t)=\left(  L_{1}(t){\small ,}L_{2}%
(t)\right)  \in\Pi_{x_{_{1}},x_{2}}$, we consider the one-dimensional payment
strategy $L(t)=L_{1}(t)+L_{2}(t)$ for the merged company. Since $\tau_{M}%
^{L}\geq\widehat{\tau}$ we conclude that $V_{M}(x_{1}+x_{2})\geq
V(\overline{x})$. Note that this is not longer true as $m>0$.
\end{remark}

\section{Properties of the optimal value function}

\label{sec:hjb}

In this section we first state some basic results concerning regularity and
growth at infinity of the optimal value function $V$ defined in (\ref{V}). We
then deduce the Hamilton-Jacobi-Bellman equation associated to this
optimization problem and see that $V$ is a viscosity solution of the HJB
equation. Moreover, we can characterize $V$ as the smallest viscosity
supersolution of this equation with a suitable growth condition. Finally, we
obtain a verification result: any viscosity supersolution which is a value
function of a family of admissible strategies is the optimal value function.

The following two lemmas are the two-dimensional counterparts of Lemmas 2.1
and 2.2 of \cite{azmuazul}.

\begin{lemma}
\label{V_GrowthCondition} For all $(x_{1},x_{2})\in\mathbf{R}_{+}^{2}$ the
optimal value function is well defined and satisfies%
\[
x_{1}+x_{2}+\frac{c_{1}+c_{2}}{q+\lambda}\leq V(x_{1},x_{2})\leq x_{1}%
+x_{2}+\frac{c_{1}+c_{2}}{q}.
\]

\end{lemma}

\textbf{Proof.} Taking strategy that pays at the beginning the whole surpluses
$x_{1}+x_{2}$ and then paying incoming premium as dividends up to the first
claim arrival $\tau_{1}$ gives the lower bound since
\[
V(x_{1},x_{2})\geq x_{1}+x_{2}+(c_{1}+c_{2})E\int_{0}^{\tau_{1}}%
e^{-qt}dt=x_{1}+x_{2}+\frac{c_{1}+c_{2}}{q+\lambda}.
\]
Similarly, observation that $L_{i}(s)\leq x_{i}+c_{i}s$ ($i=1,2$) produces the
upper bound. \vspace{3mm}\hfill$\Box$

\begin{lemma}
\label{V increasing_LocallyLip} The optimal value function $V$ is increasing,
locally Lipschitz and satisfies
\[
h\leq V(x_{1}+h,x_{2})-V(x_{1},x_{2})\leq(e^{(q+\lambda)h/c_{1}}%
-1)V(x_{1},x_{2})
\]
and%
\[
h\leq V(x_{1},x_{2}+h)-V(x_{1},x_{2})\leq(e^{(q+\lambda)h/c_{2}}%
-1)V(x_{1},x_{2})
\]
for any $(x_{1},x_{2})\in\mathbf{R}_{+}^{2}\ $and any $h>0$.
\end{lemma}

\textbf{Proof.} To prove the lower inequality in the first series of
inequalities, it suffices to consider the strategy that pays $h$ as dividends
from the surplus of the first branch and then follows the strategy
$\pi_{\epsilon}$ such that $V_{\pi_{\epsilon}}(x_{1},x_{2})\geq V(x_{1}%
,x_{2})-\epsilon$ for general $\epsilon>0$. To prove the upper inequality, it
suffices to consider the strategy that pays no dividends up the first passage
time of $x_{1}+h$ by the first surplus process $X_{1}(t)$ and the follows
strategy $\pi_{\epsilon}$ such that $V_{\pi_{\epsilon}}(x_{1}+h,x_{2})\geq
V(x_{1}+h,x_{2})-\epsilon$ for general $\epsilon>0$. The second line of
inequalities could be proved in a similar way. The proof of the fact that $V$
is increasing and locally Lipschitz follows now classical arguments.
\vspace{3mm}\hfill$\Box$

In order to obtain the Hamilton-Jacobi-Bellman (HJB) equation associated to
the optimization problem (\ref{V}), we need to state the so-called Dynamic
Programming Principle (DPP). Since $\overline{X}(t)$ is a Markov process the
proof follows the same arguments as the ones given in Lemma 1.2 of Azcue and
Muler \cite{azmu}. They use only the fact that $V$ is increasing and
continuous in $\mathbf{R}_{+}^{2}$. Obviously we need also to extend the
definition of $V$ to $\mathbf{R}^{2}$ defining $V$ as zero outside the first quadrant.

\begin{lemma}
\label{DPP}For any initial surplus $\overline{x}$ in $\mathbf{R}_{+}^{2}$ and
any stopping time $\tau$, we have%
\[%
\begin{array}
[c]{l}%
V(\overline{x})\\%
\begin{array}
[c]{ll}%
= & \sup\limits_{\overline{L}\in\Pi_{\overline{x}}}E_{\overline{x}}\left(
%TCIMACRO{\tint \limits_{0}^{\tau\wedge\widehat{\tau}}}%
%BeginExpansion
{\textstyle\int\limits_{0}^{\tau\wedge\widehat{\tau}}}
%EndExpansion
e^{-qs}dL_{1}(s)+%
%TCIMACRO{\tint \limits_{0}^{\tau\wedge\widehat{\tau}}}%
%BeginExpansion
{\textstyle\int\limits_{0}^{\tau\wedge\widehat{\tau}}}
%EndExpansion
e^{-qs}dL_{2}(s)+e^{-q(\tau\wedge\widehat{\tau})}V(X_{1}^{\overline{L}}\left(
\tau\wedge\widehat{\tau}\right)  ,X_{2}^{\overline{L}}\left(  \tau
\wedge\widehat{\tau}\right)  )\right)  .
\end{array}
\end{array}
\]

\end{lemma}

We now deduce the HJB equation assuming some regularity on $V$.

For any continuously differentiable function $u$ defined in $\mathbf{R}%
_{+}^{2}$, we define the infinitesimal generator $\widetilde{\mathcal{G}}$ of
the controlled process $\overline{X}^{\overline{L}}(t\wedge\hat{\tau})$ by
\begin{equation}%
\begin{array}
[c]{lll}%
\widetilde{\mathcal{G}} u(\overline{x}) & := & \lim_{t\rightarrow0}%
\frac{E_{\overline{x}}(e^{-q\,t}u(\overline{X}^{\overline{L}}(t\wedge\hat
{\tau}))-u(\overline{x}))}{t};
\end{array}
\label{Inf_Generator0}%
\end{equation}
see \cite[Sec. 1.4]{azmu} for details.

We will consider now the admissible strategy $\overline{L}$ in which both
branches pay dividends with constant rates $l_{1}\geq0$ and $l_{2}\geq0$
respectively until ruin time $\hat{\tau}$. Then using \cite[Rem. 1.8]{azmu} we
have
\begin{equation}%
\begin{array}
[c]{lll}%
\widetilde{\mathcal{G}}u(\overline{x}) & = & \left(  c_{1}-l_{1}\right)
u_{x_{1}}(\overline{x})+\left(  c_{2}-l_{2}\right)  u_{x_{2}}(\overline{x})\\
&  & -\left(  q+\lambda\right)  u(\overline{x})+\mathcal{I}(u)(\overline{x}),
\end{array}
\label{Inf_Generator}%
\end{equation}
where%

\begin{equation}
\mathcal{I}(u)(\overline{x})=\lambda\int_{0}^{\left(  x_{1}/b_{1}\right)
\wedge\left(  x_{2}/b_{2}\right)  }u(x_{1}-b_{1}\alpha,x_{2}-b_{2}%
\alpha)dG(\alpha) \label{Definicion_I}%
\end{equation}
and $\tau_{1}$ denotes the first claim arrival.

Assume now that $V$ is continuously differentiable. Note that from Lemma
\ref{DPP} we have:
\[
V\left(  \overline{x}\right)  \geq(l_{1}+l_{2})E_{\overline{x}}\int_{0}%
^{\tau_{1}}e^{-qt}dt+E_{\overline{x}}e^{-q\left(  t\wedge\tau_{1}\right)
}V\left(  \overline{X}^{\overline{L}}(t\wedge\tau_{1})\right)  .
\]
Thus by \eqref{Inf_Generator} dividing above inequality by $t$ and taking
$t\downarrow0$ give:
\[
\mathcal{L}(V)(\overline{x})+l_{1}\left(  1-V_{x_{1}}(\overline{x})\right)
+l_{2}\left(  1-V_{x_{2}}(\overline{x})\right)  \leq0,
\]
where%

\begin{equation}
\mathcal{L}(V)(\overline{x}):=c_{1}V_{x_{1}}(\overline{x})+c_{2}V_{x_{2}%
}(\overline{x})-\left(  q+\lambda\right)  V(\overline{x})+\mathcal{I}%
(V)(\overline{x}). \label{DefinicionL}%
\end{equation}
Taking $l_{1}=l_{2}=0$, $l_{1}\rightarrow\infty$ with $l_{2}=0$; and
$l_{2}\rightarrow\infty$ with $l_{1}=0$, we obtain%

\begin{equation}
\label{3.3b}\max\left\{  \mathcal{L}(V)(\overline{x}),1-V_{x_{1}}(\overline
{x}),1-V_{x_{2}}(\overline{x})\right\}  \leq0.\hfill
\end{equation}
We now associate to our problem the following HJB equation:%

\begin{equation}
\max\left\{  \mathcal{L}(V),1-V_{x_{1}},1-V_{x_{2}}\right\}  =0. \label{HJB}%
\end{equation}

Since the optimal value function could be not differentiable, we have to use
the notion of viscosity solution and see that $V$ is a viscosity solution of
the associated HJB equation. Let us define this notion (see for instance
Crandall and Lions \cite{cralio} and Soner \cite{soner}).

\begin{definition}
\label{Viscosity} A locally Lipschitz function $\overline{u}:\mathbf{R}%
_{+}^{2}\rightarrow\mathbf{R}$\ is a \textit{viscosity supersolution} of
(\ref{HJB})\ at $\overline{x}\in\mathbf{R}_{+}^{2}$\ if any continuously
differentiable function $\varphi:\mathbf{R}_{+}^{2}\rightarrow\mathbf{R}%
\ $with $\varphi(\overline{x})=\overline{u}(\overline{x})$ such that
$\overline{u}-\varphi$\ reaches the minimum at $\overline{x}$\ satisfies
\[
\max\left\{  \mathcal{L}(\varphi)(\overline{x}),1-\varphi_{x}(\overline
{x}),1-\varphi_{x_{2}}(\overline{x})\right\}  \leq0.\
\]
A function $\underline{u}:$ $\mathbf{R}_{+}^{2}\rightarrow\mathbf{R}$\ is a
\textit{viscosity subsolution}\ of (\ref{HJB}) at $\overline{x}\in
\mathbf{R}_{+}^{2}$\ if any continuously differentiable function
$\psi:\mathbf{R}_{+}^{2}\rightarrow\mathbf{R}\ $with $\psi(\overline
{x})=\underline{u}(\overline{x})$ such that $\underline{u}-\psi$\ reaches the
maximum at $\overline{x}$ satisfies
\[
\max\left\{  \mathcal{L}(\psi)(\overline{x}),1-\psi_{x_{1}}(\overline
{x}),1-\psi_{x_{2}}(\overline{x})\right\}  \geq0\text{.}%
\]
A function $u:\mathbf{R}_{+}^{2}\rightarrow\mathbf{R}$ which is both a
supersolution and subsolution at $\overline{x}\in\mathbf{R}_{+}^{2}$ is called
a \textit{viscosity solution }of (\ref{HJB})\ at $\overline{x}\in
\mathbf{R}_{+}^{2}$.
\end{definition}

\begin{theorem}
\label{Prop V is a viscosity supersolution} $V$ is a viscosity solution of the
HJB equation (\ref{HJB}) at any $\overline{x}\in\mathbf{R}_{+}^{2}.$
\end{theorem}

\textbf{Proof.} The proof that $V$ is a viscosity supersolution is similar to
the one in Proposition 3.1 in \cite{azmu}. We underline only crucial
adjustments that should be made in the proof.

The proof that $V$ is a viscosity supersolution follows the same arguments as
the ones used to derive (\ref{3.3b}).

The proof of the fact that $V$ is a viscosity subsolution is done by
contradiction. We will use the following convention. For any vectors
$\overline{a}=(a_{1},a_{2})$ and $\overline{b}=(b_{1},b_{2})$ we denote
$[\overline{a}, \overline{b}]:=[a_{1},b_{1}]\times[a_{2},b_{2}]$ with
$[-\infty, \overline{b}]:=[-\infty,b_{1}]\times[-\infty,b_{2}]$ and for any
$h>0$ we will write $\overline{a}\pm h:=(a_{1}\pm h, a_{2}\pm h)$. We assume
that for some fixed $x_{0}=(x_{01},x_{02})$ there exist $\epsilon>0$ and
$h\in(0,\frac{1}{2}x_{01}\wedge x_{02})$ and test function $\psi$ such that:
\[
1-\psi_{x_{i}}(\overline{x})\leq0, \quad i=1,2
\]
for $\overline{x}\in[0,\overline{x}_{0}+h]$,
\[
\mathcal{L}(\psi)(\overline{x})\leq-2q\epsilon
\]
for $\overline{x}\in[\overline{x}_{0}-h,\overline{x}_{0}+h]$,
\[
V(\overline{x})\leq\psi(\overline{x})-2\epsilon
\]
for $\overline{x}\in[-\infty,\overline{x}_{0}-h/2]\cup\{\overline{x}_{0}+h\}$.
Now the proof will go along the lines of Proposition 3.1 in \cite{azmu} by
taking $(x_{1}-x_{01})^{2}(x_{2}-x_{02})^{2}$ and $x_{01}^{2}x_{02}^{2}$
instead of $(x-x_{0})^{2}$ and $x_{0}^{2}$, respectively. One needs to modify
also the definitions $\overline{\tau}$ and $\underline{\tau}$ into the
following ones:
\[
\overline{\tau}=\inf\{t\geq0: X_{1}(t)\geq x_{01}+h\quad\text{or}\quad
X_{2}(t)\geq x_{02}+h\}
\]
and
\[
\underline{\tau}=\inf\{t\geq0: X_{1}(t)\leq x_{01}-h\quad\text{or}\quad
X_{2}(t)\leq x_{02}-h\}.
\]
In the last step (see inequality (3.21) in \cite{azmu}) we use the following
zero expectation martingale%

\begin{equation}%
\begin{array}
[c]{lll}%
\widetilde{M}_{\psi}(t) & = &
%TCIMACRO{\tsum \limits_{\substack{X_{1}(s^{-})\neq X_{1}(s)\\s\leq t}}}%
%BeginExpansion
{\textstyle\sum\limits_{\substack{X_{1}(s^{-})\neq X_{1}(s)\\s\leq t}}}
%EndExpansion
\left(  \psi(\overline{X}(s))-\psi(\overline{X}(s^{-})\right)  e^{-qs}\\
&  & -\lambda\int\limits_{0}^{t}e^{-qs}\int\limits_{0}^{\frac{X_{1}\left(
s^{-}\right)  }{b_{1}}\wedge\frac{X_{2}\left(  s^{-}\right)  }{b_{2}}}\left(
\psi(\overline{X}(s^{-})-\alpha(b_{1},b_{2}))-\psi(\overline{X}(s^{-})\right)
)dG(\alpha)ds,
\end{array}
\label{DefinicionMt}%
\end{equation}
that could be defined properly by Dynkin's formula for any test function for
viscosity subsolution $\psi$; see also \cite[Prop. 2.13]{azmu}. At the end
after some manipulations we derive that
\[
V(\overline{x}_{0})\leq\psi(\overline{x}_{0})-\epsilon
\]
which is a contradiction with the assumption that $V(\overline{x}_{0}%
)=\psi(\overline{x}_{0})$. This completes the proof. \vspace{3mm}$\Box$

\begin{remark}
\label{Muchas soluciones viscosas}The\vspace{3mm} functions $U(x_{1}%
,x_{2})=x_{1}+x_{2}+K$ are viscosity solutions of the HJB equation (\ref{HJB})
for $K>$ $(c_{1}+c_{2})/q$ because%
\[%
\begin{array}
[c]{lll}%
\mathcal{L}(U)(\overline{x}) & = & c_{1}+c_{2}-\left(  q+\lambda\right)
\left(  x_{1}+x_{2}+K\right)  +\lambda\int_{0}^{\left(  x_{1}/b_{1}\right)
\wedge\left(  x_{2}/b_{2}\right)  }(x_{1}+x_{2}-\alpha+K)dG(\alpha)\\
& \leq & c_{1}+c_{2}-(q+\lambda)\left(  x_{1}+x_{2}+K\right)  +\lambda
(x_{1}+x_{2}+K)G(\left(  x_{1}/b_{1}\right)  \wedge\left(  x_{2}/b_{2}\right)
)\\
& \leq & c_{1}+c_{2}-qK.
\end{array}
\]
So there are infinitely many viscosity solutions of the HJB equation
(\ref{HJB}).
\end{remark}

We now see that the optimal value function $V$ can be characterized as the
smallest viscosity supersolution of (\ref{HJB}) with a suitable growth
condition. We say that the function $u:\mathbf{R}_{+}^{2}\rightarrow
\mathbf{R}$ satisfies the growth condition, if there exists $K>0$ such that%

\begin{equation}
u(\overline{x})\leq K+x_{1}+x_{2}\text{ for all }\overline{x}\in\mathbf{R}%
_{+}^{2}\text{.} \tag{\bf G}\label{growth}%
\end{equation}
We have the following result.

\begin{proposition}
\label{Supersolution mayor que estrategia}Let $\overline{L}\in\Pi
_{\overline{x}}$ be any admissible strategy and let $\overline{u}$ be any
viscosity supersolution of (\ref{HJB}) satisfying growth condition
(\ref{growth}), we have that $V_{\overline{L}}(\overline{x})\leq\overline
{u}(x).$
\end{proposition}

The proof of this proposition is a straightforward extension to the
two-dimensional case of the corresponding proof of Proposition 4.4 in
\cite{azmu}, taking again the zero-expectation martingale defined in
(\ref{DefinicionMt}) and using Lemmas \ref{V_GrowthCondition}, and
\ref{V increasing_LocallyLip}.

From the last proposition and Theorem
\ref{Prop V is a viscosity supersolution}, we conclude the following
corollary.\hfill

\begin{corollary}
\label{MenorSuper} The optimal value function $V$ is the smallest viscosity
supersolution of (\ref{HJB}) satisfying growth condition (\ref{growth}).
\end{corollary}

\begin{remark}
\label{Verification_V} From Proposition
\ref{Supersolution mayor que estrategia} we can deduce the usual viscosity
verification result: If the value function $V_{\pi}$ for some strategy $\pi$
or a limit of value functions $\lim_{n\rightarrow\infty}V_{\pi_{n}}$ for some
strategies $\pi_{n}$ is a viscosity supersolution of (\ref{HJB}), then it is
the optimal value function (\ref{V}). Note that by definition $V$ is a limit
of value functions. However, we expect $V$ to be a value function of a family
of admissible strategies with some particular structure. This problem will be
analyzed in next sections.
\end{remark}

\section{Optimal strategy}

\label{sec: optimalstrategy}

First in this section, we introduce some special families of admissible
strategies which depends only on the current surplus levels.

Assuming that $V$ is differentiable in $\mathbf{R}_{+}^{2}$, the equation
(\ref{HJB}) suggests how the dividends are paid depending on the current
surplus $\overline{x}=(x_{1},x_{2})\in\mathbf{R}_{+}^{2}$:

\begin{itemize}
\item If the current surplus is in the set%
\[
\mathcal{C}^{\ast}=\left\{  \overline{x}\in\mathbf{R}_{+}^{2}:\mathcal{L}%
(V)(\overline{x})=0\text{, }1-V_{x_{1}}(\overline{x})<0\text{, }1-V_{x_{2}%
}(\overline{x})<0\text{ }\right\}
\]
no dividends are paid.

\item If the current surplus is in the set%
\[
\mathcal{B}_{0}^{\ast}=\left\{  \overline{x}\in\mathbf{R}_{+}^{2}%
:\mathcal{L}(V)(\overline{x})<0\text{, }1-V_{x_{1}}(\overline{x})=0\text{,
}1-V_{x_{2}}(\overline{x})=0\text{ }\right\}
\]
both branches pay a lump sum as dividends.

\item If the current surplus is in the set%
\[
\mathcal{B}_{1}^{\ast}=\left\{  \overline{x}\in\mathbf{R}_{+}^{2}%
:\mathcal{L}(V)(\overline{x})<0\text{, }1-V_{x_{1}}(\overline{x})=0\text{,
}1-V_{x_{2}}(\overline{x})<0\text{ }\right\}
\]
the first branch pays a lump sum as dividends.

\item If the current surplus is in the set%
\[
\mathcal{B}_{2}^{\ast}=\left\{  \overline{x}\in\mathbf{R}_{+}^{2}%
:\mathcal{L}(V)(\overline{x})<0\text{, }1-V_{x_{1}}(\overline{x})<0\text{,
}1-V_{x_{2}}(\overline{x})=0\text{ }\right\}
\]
the second branch pays a lump sum as dividends.

\item If the current surplus is in the set%
\[
\mathcal{A}_{0}^{\ast}=\left\{  \overline{x}\in\mathbf{R}_{+}^{2}%
:\mathcal{L}(V)(\overline{x})=0\text{, }1-V_{x_{1}}(\overline{x})=0\text{,
}1-V_{x_{2}}(\overline{x})=0\text{ }\right\}
\]
both branches pay their incoming premiums as dividends.

\item If we define $\mathcal{A}_{1}^{\ast}$ as the boundary between
$\mathcal{B}_{1}^{\ast}$ and $\mathcal{C}^{\ast}$, we would have%
\[
\mathcal{A}_{1}^{\ast}=\left\{  \overline{x}\in\mathbf{R}_{+}^{2}%
:\mathcal{L}(V)(\overline{x})=0\text{, }1-V_{x_{1}}(\overline{x})=0\text{,
}1-V_{x_{2}}(\overline{x})<0\text{ }\right\}  .
\]
If the current surplus is in the set $\mathcal{A}_{1}^{\ast}$, the first
branch could pay some part of its incoming premium as dividends.

\item Analogously, if we define $\mathcal{A}_{2}^{\ast}$ as the boundary
between $\mathcal{B}_{2}^{\ast}$ and $\mathcal{C}^{\ast}$, we would have%
\[
\mathcal{A}_{2}^{\ast}=\left\{  \overline{x}\in\mathbf{R}_{+}^{2}%
:\mathcal{L}(V)(\overline{x})=0\text{, }1-V_{x_{1}}(\overline{x})<0\text{,
}1-V_{x_{2}}(\overline{x})=0\text{ }\right\}  \text{.}%
\]
If the current surplus is in the set $\mathcal{A}_{2}^{\ast}$, the second
branch could pay some part of its incoming premium as dividends.
\end{itemize}

If $V$ is not differentiable at some points, still these sets could be defined
in a viscosity sense.

Let us define $\mathcal{M}$ as the half line in $\mathbf{R}_{+}^{2}$ that
contains the origin with slope $b_{2}/b_{1}$. Note that, when the current
surplus is outside $\mathcal{M}$, there exists the possibility that the ruin
occur leaving one branch with positive surplus; if the current surplus is in
the set%

\[
D^{1}:=\left\{  \overline{x}\in\mathbf{R}_{+}^{2}:\left(  b_{2}/b_{1}\right)
x_{1}>x_{2}\right\}  ,
\]
then the first branch will be the one with eventual positive surplus at ruin
time and conversely, if the current surplus is in the set%

\[
D^{2}:=\left\{  \overline{x}\in\mathbf{R}_{+}^{2}:\left(  b_{2}/b_{1}\right)
x_{1}<x_{2}\right\}  ,
\]
the same could happen with the second branch.

\subsection{Optimal strategy in $D^{1}$\label{Optimal strategy in D1}}

Given a initial surplus $(x_{1},x_{2})\in\mathbf{R}_{+}^{2},$ we define the
set $\widehat{\Pi}_{x_{1},x_{2}}\subset$ $\Pi_{x_{1},x_{2}}$ in the following
way: If $(x_{1},x_{2})\in$ $D^{2}\cup\mathcal{M}$ then $\widehat{\Pi}%
_{x_{1},x_{2}}$ is the set of all the admissible strategies in which the
controlled process remains in the set $D^{2}\cup\mathcal{M}$ until ruin time;
if $(x_{1},x_{2})\in$ $D^{1},$ then $\widehat{\Pi}_{x_{1},x_{2}}$ is the set
of all the admissible strategies in which the first branch pays immediately at
least $x_{1}-\left(  b_{1}/b_{2}\right)  x_{2}$ as dividends and afterwards
the controlled process remains in the set $D^{2}\cup\mathcal{M}$. From
(\ref{X_Controlado}), we get%

\begin{equation}
\widehat{\Pi}_{x_{1},x_{2}}=\left\{  \left(  L_{1}(t),L_{2}(t)\right)  \in
\Pi_{x_{1},x_{2}}:L_{1}(t)\geq x_{1}-\tfrac{b_{1}}{b_{2}}x_{2}+(c_{1}%
-\tfrac{b_{1}}{b_{2}}c_{2})t+\tfrac{b_{1}}{b_{2}}L_{2}(t)\text{ for }%
t\geq0\right\}  . \label{PiBarra}%
\end{equation}

\begin{proposition}
\label{ProposicionEstrategias} We have that $V(x_{1},x_{2})=\sup_{\overline
{L}\in\widehat{\Pi}_{x_{1},x_{2}}}V_{\overline{L}}(x_{1},x_{2}).$
\end{proposition}

\textbf{Proof.} Given any admissible $\overline{L}^{0}=(L_{1}^{0},L_{2}%
^{0})\in\Pi_{x_{1},x_{2}}$, let us define
\begin{equation}
\overline{L}^{1}(t)=\left(  \max\{x_{1}-\tfrac{b_{1}}{b_{2}}x_{2}%
+(c_{1}-\tfrac{b_{1}}{b_{2}}c_{2})t+\tfrac{b_{1}}{b_{2}}L_{2}^{0}(t),L_{1}%
^{0}(t)\},L_{2}^{0}(t)\right)  . \label{DefinicionL1}%
\end{equation}
We have that $\overline{L}^{1}\in\widehat{\Pi}_{x_{1},x_{2}}$ because
$L_{1}^{1}(t)\ $and $L_{2}^{1}(t)$ are predictable, positive and increasing
and also
\[
L_{1}^{1}(t)\geq x_{1}-\tfrac{b_{1}}{b_{2}}x_{2}+(c_{1}-\tfrac{b_{1}}{b_{2}%
}c_{2})t+\tfrac{b_{1}}{b_{2}}L_{2}^{1}(t).
\]

Let us prove now that the ruin time $\widehat{\tau}^{\overline{L}^{0}}$
corresponding to the admissible strategy $\overline{L}^{0}$ is less than or
equal to $\widehat{\tau}^{\overline{L}^{1}}$. By definition, we have that%
\[
x_{1}+c_{1}t-%
%TCIMACRO{\dsum \limits_{i=1}^{N_{t}}}%
%BeginExpansion
{\displaystyle\sum\limits_{i=1}^{N_{t}}}
%EndExpansion
b_{1}U_{i}-L_{1}^{0}(t)\geq0
\]

and
\[
x_{2}+c_{2}t-%
%TCIMACRO{\dsum \limits_{i=1}^{N_{t}}}%
%BeginExpansion
{\displaystyle\sum\limits_{i=1}^{N_{t}}}
%EndExpansion
b_{2}U_{i}-L_{2}^{0}(t)\geq0
\]
for any $t\leq\widehat{\tau}^{\overline{L}^{0}}$. Thus $X_{2}^{\overline
{L}^{1}}(t)=X_{2}^{\overline{L}^{0}}(t)\geq0,$ and
\[%
\begin{array}
[c]{lll}%
X_{1}^{\overline{L}^{1}}(t) & = & x_{1}+c_{1}t-%
%TCIMACRO{\dsum \limits_{i=1}^{N_{t}}}%
%BeginExpansion
{\displaystyle\sum\limits_{i=1}^{N_{t}}}
%EndExpansion
b_{1}U_{i}-L_{1}^{1}(t)\\
& = & \left\{
\begin{array}
[c]{ccc}%
X_{1}^{\overline{L}^{0}}(t) & \text{if } & x_{1}-\tfrac{b_{1}}{b_{2}}%
x_{2}+(c_{1}-\tfrac{b_{1}}{b_{2}}c_{2})t+\tfrac{b_{1}}{b_{2}}L_{2}^{0}(t)\leq
L_{1}^{0}(t)\\
\tfrac{b_{1}}{b_{2}}X_{2}^{\overline{L}^{0}}(t) & \text{if} & x_{1}%
-\tfrac{b_{1}}{b_{2}}x_{2}+(c_{1}-\tfrac{b_{1}}{b_{2}}c_{2})t+\tfrac{b_{1}%
}{b_{2}}L_{2}^{0}(t)>L_{1}^{0}(t)
\end{array}
\right. \\
& \geq & 0
\end{array}
\]

for any $t\leq\widehat{\tau}^{\overline{L}^{0}}$.

\bigskip From (\ref{DefinicionL1}) we have that $L_{1}^{1}(t)+L_{2}^{1}(t)\geq
L_{1}^{0}(t)+L_{2}^{0}(t)$ for $0\leq t<$ $\widehat{\tau}^{\overline{L}^{1}}$,
so we obtain
\[
V_{\overline{L}^{0}}(\overline{x})\leq V_{\overline{L}^{1}}(\overline{x})
\]
and the result follows. \vspace{3mm}\hfill$\Box$

\begin{corollary}
\label{corolario lineal en D1}If $(x_{1},x_{2})\in$ $D^{1},$ then
$V(x_{1},x_{2})=x_{1}-\left(  b_{1}/b_{2}\right)  x_{2}+V(\left(  b_{1}%
/b_{2}\right)  x_{2},x_{2})$.
\end{corollary}

\textbf{Proof.} By (\ref{PiBarra}), we have that $\overline{L}^{0}=\left(
L_{1}(t),L_{2}(t)\right)  \in\widehat{\Pi}_{x_{1},x_{2}}$ if and only if
\[
\overline{L}^{1}=\left(  L_{1}(t)-x_{1}-\left(  b_{1}/b_{2}\right)
x_{2},L_{2}(t)\right)  \in\widehat{\Pi}_{\left(  b_{1}/b_{2}\right)
x_{2},x_{2}}.
\]
Since we have that $V_{\overline{L}^{0}}(x_{1},x_{2})=x_{1}-\left(
b_{1}/b_{2}\right)  x_{2}+V_{\overline{L}^{1}}(\left(  b_{1}/b_{2}\right)
x_{2},x_{2}),$ the result follows by Proposition \ref{ProposicionEstrategias}.
\vspace{3mm}\hfill$\Box$

From the last Corollary we can deduce that if $\overline{x}\in D^{1}$ then the
optimal value function is a limit of strategies in which the first branch
immediately pays $\left(  x_{1}-\left(  b_{1}/b_{2}\right)  x_{2}\right)  $ as
dividends, so the current surplus goes immediately to $\mathcal{M}$ in the
horizontal direction, and afterwards the controlled process remains in
$D^{2}\cup\mathcal{M}$ until ruin time. We can ask ourselves if an analogous
result holds for $\overline{x}\in D^{2}$, that is whether the optimal value
function is a limit of strategies in which the second branch immediately pays
$x_{2}-\left(  b_{2}/b_{1}\right)  x_{1}$ as dividends -so the current surplus
goes immediately to $\mathcal{M}$ in the vertical direction- and afterwards
the controlled process remains in $\mathcal{M}$ until ruin time. We will see
in the next section that in general this holds true only in the case in which
$c_{2}/c_{1}=b_{2}/b_{1}$.

\subsection{$\mathcal{M}$-strategies}

Given an initial surplus $(x_{1},x_{2})\in\mathbf{R}_{+}^{2}$, we define the
set $\widetilde{\Pi}_{x_{1},x_{2}}\subset\widehat{\Pi}_{x_{1},x_{2}}$
$\subset\Pi_{x_{1},x_{2}}$ as the set all the admissible strategies which pay
immediately dividends in the following way: the first branch pays immediately
$x_{1}-\left(  b_{1}/b_{2}\right)  x_{2}$ as dividends if $(x_{1},x_{2})\in$
$D^{1}\cup\mathcal{M}\ $or the second branch pays immediately $x_{2}-\left(
b_{2}/b_{1}\right)  x_{1}$ as dividends if $(x_{1},x_{2})\in$ $D^{2}$;
afterwards the controlled process remains in the set $\mathcal{M}$ until ruin time.

Let us define
\begin{equation}
\widetilde{V}(x_{1},x_{2}):=\sup_{\overline{L}\in\widetilde{\Pi}_{x_{1},x_{2}%
}}V_{\overline{L}}(x_{1},x_{2}). \label{VRulo}%
\end{equation}
Let us also define $\overline{V}$ as the value function of the best strategy
with initial surplus $(\left(  b_{1}/b_{2}\right)  x_{2},x_{2})\in\mathcal{M}$
whose controlled process remains in the set $\mathcal{M}$ until ruin time,
that is
\begin{equation}
\overline{V}(x_{2}):=\widetilde{V}(\tfrac{b_{1}}{b_{2}}x_{2},x_{2}).
\label{AuxiliaryOptimizationProblem}%
\end{equation}
By definition, we have that%
\[
\widetilde{V}(x_{1},x_{2})=\left(  x_{1}-\tfrac{b_{1}}{b_{2}}x_{2}%
+\overline{V}(x_{2})\right)  I_{D^{1}\cup\mathcal{M}}(x_{1},x_{2})+\left(
x_{2}-\tfrac{b_{2}}{b_{1}}x_{1}+\overline{V}(\tfrac{b_{2}}{b_{1}}%
x_{1})\right)  I_{D^{2}}(x_{1},x_{2}).
\]

\hfill

In order to find $\overline{V}$, let us consider the following auxiliary
one-dimensional optimization problem. Let $\Pi_{x_{2}}$ be the set of
admissible dividend payment strategy corresponding to the one-dimensional
compound Poisson process with initial surplus $x_{2}\geq0$, slope $c_{2}$,
claim arrivals given by the Poisson process $N_{t}$ and claim-sizes
$b_{2}U_{i}$. Given any $L_{2}\in\Pi_{x_{2}},$ we define%

\[
W_{L_{2}}(x_{2}):=E_{x_{2}}\left(  \int_{0}^{\widehat{\tau}}e^{-qs}%
(c_{1}-\tfrac{b_{1}}{b_{2}}c_{2})ds+(1+\tfrac{b_{1}}{b_{2}})\int_{0}%
^{\widehat{\tau}}e^{-qs}dL_{s}^{2}\right)  ,
\]
and%

\begin{equation}
\overline{W}(x_{2})=\sup_{L_{2}\in\Pi_{x_{2}}}W_{L_{2}}(x_{2}).
\label{Definicion Wbarra}%
\end{equation}

Note that, since any $\overline{L}=(L_{1},L_{2})\in\widetilde{\Pi}%
_{(b_{1}/b_{2})x_{2},x_{2}}\ $ satisfies
\[
L_{1}(t)=\tfrac{b_{1}}{b_{2}}L_{2}(t)+(c_{1}-\tfrac{b_{1}}{b_{2}}c_{2})t,
\]
then we have that
\[
V_{\overline{L}}(\tfrac{b_{1}}{b_{2}}x_{2},x_{2})=W_{L_{2}}(x_{2}).
\]
Thus $\overline{V}$ and $\overline{W}$ coincide which is stated formally in
next proposition.

\begin{proposition}
\label{Optimizacion_Problema_Auxiliar} For any $x_{2}\geq0$ we have
$\overline{V}(x_{2})=\overline{W}(x_{2})$.
\end{proposition}

We study now the optimization problem (\ref{Definicion Wbarra}). The function
$\overline{W}$ can be seen as the optimal value function (up to a constant) of
the dividend problem with a reward for avoiding early ruin; see for instance
Thonhauser and Albrecher \cite{ThonAlb}. In effect, we can write%
\[
\frac{\overline{W}(x_{2})}{1+\tfrac{b_{1}}{b_{2}}}=\sup_{L_{2}\in\Pi_{x_{2}}%
}E_{x_{2}}\left(  \int_{0}^{\tau}e^{-qs}\kappa ds+\int_{0}^{\tau}e^{-qs}%
dL_{s}^{2}\right)  \text{,}%
\]
where%
\[
\kappa=\frac{c_{1}-\left(  b_{1}/b_{2}\right)  c_{2}}{1+b_{1}/b_{2}}.
\]
From \cite{ThonAlb}, we have that $\overline{W}$ has the following associated
HJB equation,%
\begin{equation}
\max\{1+\tfrac{b_{1}}{b_{2}}-\overline{W}^{\prime}(x_{2}),\overline
{\mathcal{L}}(\overline{W})(x_{2})\}=0, \label{HJB_Auxiliar}%
\end{equation}
where%

\[
\overline{\mathcal{L}}(\overline{W})(x_{2}):=c_{2}\overline{W}^{\prime}%
(x_{2})-(\lambda+q)\overline{W}(x_{2})+\lambda%
%TCIMACRO{\tint \limits_{0}^{x_{2}/b_{2}}}%
%BeginExpansion
{\textstyle\int\limits_{0}^{x_{2}/b_{2}}}
%EndExpansion
\overline{W}(x_{2}-b_{2}\alpha)dG(\alpha)+c_{1}-\tfrac{b_{1}}{b_{2}}c_{2}.
\]

We can state the following proposition for the optimal one-dimensional reward
problem. The proof is similar to the one-dimensional optimization dividend
problem without reward (see for instance Azcue and Muler \cite{azmu},
Propositions 3.1 and 4.4).

\begin{proposition}
$\overline{W}$ is a locally Lipschitz viscosity solution of
(\ref{HJB_Auxiliar}). Moreover, it is the smallest viscosity solution of this
equation for all $x_{2}>0$ satisfying the growth condition $\overline{W}%
(x_{2})\leq K+\left(  1+b_{1}/b_{2}\right)  x_{2}$ for some $K>0$.
\end{proposition}

From the previous proposition, we can deduce the following verification
result: Given a family of admissible strategies $\overline{\pi}=\left\{
L_{2}^{x_{2}}\in\Pi_{x_{2}}\text{ for any }x_{2}\geq0\right\}  $, we define
the \textit{value function} $W_{\overline{\pi}}:\mathbf{R}_{+}\rightarrow
\mathbf{R}_{+}$ as $W_{\overline{\pi}}(x_{2})=W_{L_{2}^{x_{2}}}(x_{2})$. If
the function $W_{\overline{\pi}}(x_{2})$ is a viscosity supersolution of
(\ref{HJB_Auxiliar}) for each $x_{2}>0$, then $W_{\overline{\pi}}$ coincides
with $\overline{W}$.

We also have the following result about the optimal dividend strategy of the
problem (\ref{Definicion Wbarra}), the proof is similar to the case without
reward, see for instance Azcue and Muler \cite{azmu}, Sections 5.1 and 5.2.

\begin{proposition}
\label{BandPartitionLine}The dividend strategy that optimizes
(\ref{Definicion Wbarra}) is a band dividend strategy depending on the sets
$\overline{\mathcal{A}}$, $\overline{\mathcal{B}}$ and $\overline{\mathcal{C}%
}$, where
\[
\overline{\mathcal{A}}:\mathcal{=}\left\{  x_{2}\geq0:\left(  c_{1}%
+c_{2}\right)  -(\lambda+q)\overline{W}(x_{2})+\lambda%
%TCIMACRO{\tint _{0}^{x_{2}/b_{2}}}%
%BeginExpansion
{\textstyle\int_{0}^{x_{2}/b_{2}}}
%EndExpansion
\overline{W}(x_{2}-b_{2}\alpha)dG(\alpha)=0\right\}  \text{ }%
\]
is closed and bounded,
\[
\overline{\mathcal{B}}:=\left\{  x_{2}\geq0:\overline{W}^{\prime}%
(x_{2})=1+\tfrac{b_{1}}{b_{2}},\left(  c_{1}+c_{2}\right)  -(\lambda
+q)\overline{W}(x_{2})+\lambda%
%TCIMACRO{\tint \limits_{0}^{x_{2}/b_{2}}}%
%BeginExpansion
{\textstyle\int\limits_{0}^{x_{2}/b_{2}}}
%EndExpansion
\overline{W}(x_{2}-b_{2}\alpha)dG(\alpha)<0\right\}
\]
is left-open,%
\[
\overline{\mathcal{C}}=\mathbf{R}_{+}-\overline{\mathcal{A}}\cup
\overline{\mathcal{B}}%
\]
is right-open, and there exists $\widetilde{x}_{2}$ such that $(\widetilde
{x}_{2},\infty)\subset\overline{\mathcal{B}}$. We have the following cases:

\begin{itemize}
\item if the current surplus is in $\overline{\mathcal{A}}$ the incoming
premium $c_{2}$ is paid as dividends;

\item if the current surplus is in $\overline{\mathcal{B}}$ a positive amount
of money is paid as dividends in order to bring the surplus process back to
$\overline{\mathcal{A}}$;

\item if the surplus is in $\overline{\mathcal{C}}$, no dividends are paid.
\end{itemize}
\end{proposition}

Let us go back to problem (\ref{VRulo}) in $\mathcal{M}$ . Using that
\[
\widetilde{V}(\tfrac{b_{1}}{b_{2}}x_{2},x_{2})=\overline{V}(x_{2}%
)=\overline{W}(x_{2}),
\]
we can split the half line $\mathcal{M}$ into three sets:%
\[
\mathcal{A}_{\mathcal{M}}:=\left\{  \left(  \tfrac{b_{1}}{b_{2}}x_{2}%
,x_{2}\right)  :x_{2}\in\overline{\mathcal{A}}\right\}  ,
\]%
\[
\mathcal{B}_{\mathcal{M}}:=\left\{  \left(  \tfrac{b_{1}}{b_{2}}x_{2}%
,x_{2}\right)  :x_{2}\in\overline{\mathcal{B}}\right\}  ,
\]
and
\begin{equation}
\mathcal{C}_{\mathcal{M}}:=\left\{  \left(  \tfrac{b_{1}}{b_{2}}x_{2}%
,x_{2}\right)  :x_{2}\in\overline{\mathcal{C}}\right\}  . \label{CM}%
\end{equation}

From Proposition \ref{BandPartitionLine} we conclude that the optimal dividend
strategy for the auxiliary problem (\ref{VRulo}) with initial surplus in
$\mathcal{M}$ is the following:

\begin{itemize}
\item if $\left(  \left(  b_{1}/b_{2}\right)  x_{2},x_{2}\right)  \in$
$\mathcal{A}_{\mathcal{M}}$ both branches pay the incoming premium as
dividends until next claim,

\item if $\left(  \left(  b_{1}/b_{2}\right)  x_{2},x_{2}\right)
\in\mathcal{B}_{\mathcal{M}}$, the second branch pays a positive amount $m$ of
money and the first branch pays $\left(  b_{1}/b_{2}\right)  m,$ where $m$ is
the minimal amount that brings the surplus process back to $\mathcal{A}%
_{\mathcal{M}}$,

\item and finally if $\left(  \left(  b_{1}/b_{2}\right)  x_{2},x_{2}\right)
\in\mathcal{C}_{\mathcal{M}}$, the second branch pays no dividends and the
first branch pays $(c_{1}-\left(  b_{1}/b_{2}\right)  c_{2})$ as dividend rates.
\end{itemize}

\begin{remark}
\label{Ir a la Diagonal}\bigskip If the optimal value function $V$ coincides
with the function $\widetilde{V}\ $defined in (\ref{VRulo}), the optimal
strategy would be given by
\[
\mathcal{B}_{1}^{\ast}=\left\{  \overline{x}\in D^{1}:x_{2}\notin
\mathcal{B}_{\mathcal{M}}\right\}  \text{, }\mathcal{B}_{2}^{\ast}=\left\{
\overline{x}\in D^{2}:\left(  b_{2}/b_{1}\right)  x_{1}\notin\mathcal{B}%
_{\mathcal{M}}\right\}  ,
\]%
\[
\mathcal{B}_{0}^{\ast}=\left(  D^{1}-\mathcal{B}_{1}^{\ast}\right)
\cup\mathcal{B}_{\mathcal{M}}\cup\left(  D^{2}-\mathcal{B}_{2}^{\ast}\right)
,
\]
$\mathcal{A}_{0}^{\ast}=\mathcal{A}_{\mathcal{M}}$, $\mathcal{A}_{1}^{\ast
}=\mathcal{C}_{\mathcal{M}}$ and $\mathcal{A}_{2}^{\ast}=$ $\mathcal{C}^{\ast
}=\varnothing$. Moreover,

\begin{itemize}
\item in $\mathcal{A}_{0}^{\ast}$ both the incoming premiums are paid as dividends;

\item in $\mathcal{A}_{1}^{\ast}$ the first branch pays dividends at rate
$c_{1}-\left(  b_{1}/b_{2}\right)  c_{2}$ (and so the surplus process remains
in $\mathcal{M}$);

\item in $D^{1}$ the first branch pays immediately $\left(  x_{1}-\left(
b_{1}/b_{2}\right)  x_{2}\right)  $, so the current surplus goes immediately
to $\mathcal{M}$ in the horizontal direction;

\item in $D^{2}$ the second branch pays immediately $\left(  x_{2}-\left(
b_{2}/b_{1}\right)  x_{1}\right)  $, so the current surplus goes immediately
to $\mathcal{M}$ in the vertical direction.
\end{itemize}

In the particular case that $c_{1}=\left(  b_{1}/b_{2}\right)  c_{2}$, the
first branch does not need to pay dividends in $\mathcal{C}_{\mathcal{M}}$ to
remain in $\mathcal{M}$, so $\mathcal{C}^{\ast}=\mathcal{A}_{1}^{\ast
}=\mathcal{C}_{\mathcal{M}}.$
\end{remark}

\subsection{$\widetilde{V}$ is the optimal value function in the case
$c_{2}/c_{1}=b_{2}/b_{1}$\label{optimalStrategyEquality}}

In the next propositions we prove that the optimal value function $V$
coincides with the function $\widetilde{V}\ $defined in (\ref{VRulo}).

\begin{proposition}
\label{ProblemaAuxiliar} For any $x_{2}\geq0$ we have that $V(\left(
b_{1}/b_{2}\right)  x_{2},x_{2})=\overline{V}(x_{2})=\overline{W}(x_{2}).$
\end{proposition}

\textbf{Proof.} Given any admissible $\overline{L}^{0}=(L_{1}^{0},L_{2}%
^{0})\in\Pi_{\left(  b_{1}/b_{2}\right)  x_{2},x_{2}}$, let us define
$\overline{L}^{1}\in\Pi_{\left(  b_{1}/b_{2}\right)  x_{2},x_{2}}$ as
\[
\overline{L}^{1}(t)=\left(  L_{1}^{0}(t)+L_{2}^{0}(t)\right)  \left(
b_{1},b_{2}\right)
\]
Note that $L_{1}^{1}(t)+L_{2}^{1}(t)=L_{1}^{0}(t)+L_{2}^{0}(t)$ and
$\overline{X}^{\overline{L}^{1}}(t)\in\mathcal{M}$ for $t<$ $\widehat{\tau
}^{\overline{L}^{1}}$, so $\overline{L}^{1}\in\widetilde{\Pi}_{\left(
b_{1}/b_{2}\right)  x_{2},x_{2}}$ . It is easy to check that the ruin time
$\widehat{\tau}^{\overline{L}^{0}}$ corresponding to the admissible strategy
$\overline{L}^{0}$ is less than or equal to $\widehat{\tau}^{\overline{L}^{1}%
}$. Therefore,
\[
V_{\overline{L}^{0}}(\tfrac{b_{1}}{b_{2}}x_{2},x_{2})\leq V_{\overline{L}^{1}%
}(\tfrac{b_{1}}{b_{2}}x_{2},x_{2})
\]
and so the result follows. \vspace{3mm}$\Box$

\begin{proposition}
\label{V_enD1} The optimal value function $V=\widetilde{V}.$
\end{proposition}

\textbf{Proof. } By Proposition \ref{ProblemaAuxiliar} and
(\ref{Optimizacion_Problema_Auxiliar}) we have that $V=\widetilde{V}$ in
$\mathcal{M}$. Then, by Proposition \ref{corolario lineal en D1} we obtain
that $V=\widetilde{V}$ in $D^{1}$; the result in $D^{2}\ $follows by symmetry.
\vspace{3mm}$\Box$

As a consequence of the previous proposition, the optimal strategy is the one
described in Remark \ref{Ir a la Diagonal}. As an example, consider the
claim-size distribution gamma
\[
G(x)=1-(1+\frac{x}{2})e^{-\frac{x}{2}}\text{,}%
\]
the parameters $b_{1}=b_{2}=0.5$, $c_{1}=c_{2}=21.4$, $q=0.1$ and $\lambda
=10$. In Section 6.2.1 of \cite{azmu}, it is shown that $\overline
{\mathcal{A}}=\left\{  0,10.22\right\}  $, $\overline{\mathcal{B}%
}=(0,1.803]\cup(10.22,\infty)$ and $\overline{\mathcal{C}}=(1.803,10.22)$, the
optimal strategy is depicted in Figure 1.%

\[%
%TCIMACRO{\FRAME{itbpFU}{2.463in}{2.4716in}{0in}{\Qcb{Figure 1.}}%
%{}{figure1.eps}{\special{ language "Scientific Word";  type "GRAPHIC";
%width 2.463in;  height 2.4716in;  depth 0in;  original-width 4.3284in;
%original-height 4.3439in;  cropleft "0";  croptop "1";  cropright "1";
%cropbottom "0";  filename 'Figure1.eps';file-properties "XNPEU";}}}%
%BeginExpansion
{\parbox[b]{2.463in}{\begin{center}
\includegraphics[
height=2.4716in,
width=2.463in
]%
{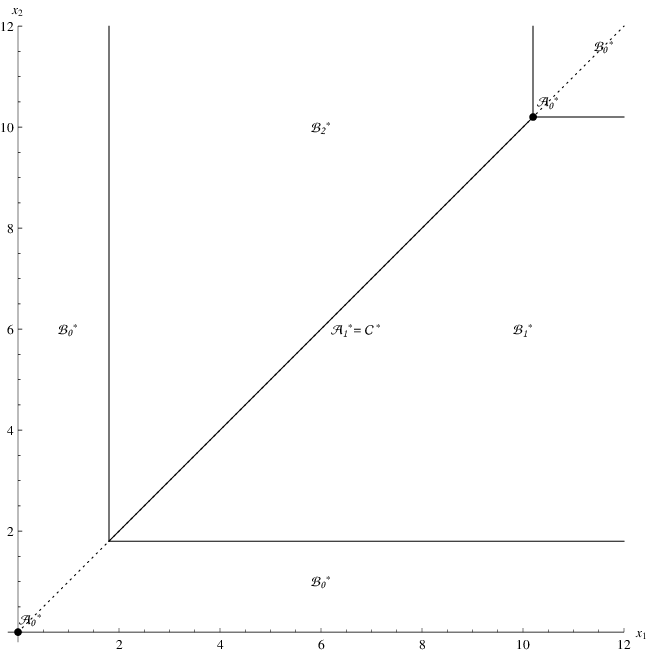}%
\\
Figure 1.
\end{center}}}%
%EndExpansion
\]

\begin{remark}
$V$ is not differentiable in $C_{\mathcal{M}}$ because if it were
differentiable at a point $\left(  \left(  b_{1}/b_{2}\right)  x_{2}%
,x_{2}\right)  \in$ $C_{\mathcal{M}}$ , then both $V_{x_{1}}\left(  \left(
b_{1}/b_{2}\right)  x_{2},x_{2}\right)  $ and $V_{x_{2}}\left(  \left(
b_{1}/b_{2}\right)  x_{2},x_{2}\right)  $ should be one; then
\[
\overline{W}%
%TCIMACRO{\U{b4}}%
%BeginExpansion
\acute{}%
%EndExpansion
(x_{2})=1+\tfrac{b_{1}}{b_{2}}%
\]
that implies $x_{2}\notin\overline{\mathcal{C}}$ which is a contradiction.
\end{remark}

\bigskip

\subsection{$\widetilde{V}$ is not the optimal value function in the case
$c_{2}/c_{1}<b_{2}/b_{1}$\label{optimalStrategyIneEquality}}

In the next proposition we show, that except in a very special case, the
function $\widetilde{V}\ $defined in (\ref{VRulo}) is never the optimal value
function of the optimization problem (\ref{V}).

\begin{proposition}
\label{Vrulo_NeverOptimal}In the case that $\mathcal{C}_{\mathcal{M}}$ defined
in (\ref{CM}) is not empty, the function $\widetilde{V}$ is \textbf{not} a
viscosity solution of (\ref{HJB}) at all the points in the first quadrant.
\end{proposition}

\textbf{Proof.} Let us take any point $\left(  x_{1}^{0},x_{2}^{0}\right)
\in\mathcal{C}_{\mathcal{M}}$ with $x_{2}^{0}=\left(  b_{2}/b_{1}\right)
x_{1}^{0}$ such that $\overline{W}$ is differentiable at $x_{2}^{0}$ (since
$\overline{W}$ is locally Lipschitz and $\overline{\mathcal{C}}$ is
right-open, the set of points in $\overline{\mathcal{C}}$ where $\overline
{W}^{\prime}$ exists has full measure). By definition of $\overline
{\mathcal{C}}$,
\[
\overline{W}^{\prime}(x_{2}^{0})-1-\tfrac{b_{1}}{b_{2}}>0,
\]
and%

\[
\overline{\mathcal{L}}(\overline{W})(x_{2}^{0})=c_{2}\overline{W}^{\prime
}(x_{2}^{0})-(\lambda+q)\overline{W}(x_{2}^{0})+\lambda%
%TCIMACRO{\tint \limits_{0}^{x_{2}^{0}/b_{2}}}%
%BeginExpansion
{\textstyle\int\limits_{0}^{x_{2}^{0}/b_{2}}}
%EndExpansion
\overline{W}(x_{2}^{0}-b_{2}\alpha)dG(\alpha)+c_{1}-\tfrac{b_{1}}{b_{2}}%
c_{2}=0.
\]
This implies that for any $x_{2}>x_{2}^{0}$ (and so $\left(  x_{1}^{0}%
,x_{2}\right)  \in D^{2}$) we have,
\[
\widetilde{V}\left(  x_{1}^{0},x_{2}\right)  =x_{2}-x_{2}^{0}+\overline
{W}(x_{2}^{0})=x_{2}-\tfrac{b_{2}}{b_{1}}x_{1}^{0}+\overline{W}(\frac{b_{2}%
}{b_{1}}x_{1}^{0})
\]
and so it is differentiable at all these points with%
\[
\widetilde{V}_{x_{1}}\left(  x_{1}^{0},x_{2}\right)  =-\tfrac{b_{2}}{b_{1}%
}+\tfrac{b_{2}}{b_{1}}\overline{W}^{\prime}(x_{2}^{0})>1\text{ and }%
\widetilde{V}_{x_{2}}\left(  x_{1}^{0},x_{2}\right)  =1.
\]
\newline Then, since $b_{1}c_{2}<b_{2}c_{1}$ and $x_{2}/b_{2}>x_{2}^{0}%
/b_{2}=x_{1}^{0}/b_{1}$,
\[%
\begin{array}
[c]{lll}%
\mathcal{L}(\widetilde{V})(x_{1}^{0},x_{2}) & = & c_{1}(-\tfrac{b_{2}}{b_{1}%
}+\tfrac{b_{2}}{b_{1}}\overline{W}^{\prime}(x_{2}^{0}))+c_{2}-\left(
q+\lambda\right)  (x_{2}-x_{2}^{0}+\overline{W}(x_{2}^{0}))\\
&  & +\lambda\int_{0}^{x_{2}^{0}/b_{2}}(x_{2}-x_{2}^{0}+\overline{W}(x_{2}%
^{0}-b_{2}\alpha))dG(\alpha)\\
& = & \overline{\mathcal{L}}(\overline{W})(x_{2}^{0})+(\tfrac{c_{1}b_{2}%
}{b_{1}}-c_{2})\left(  \overline{W}^{\prime}(x_{2}^{0})-1-\tfrac{b_{1}}{b_{2}%
}\right)  -(x_{2}-x_{2}^{0})(q+\lambda)\\
&  & +\lambda(x_{2}-x_{2}^{0})G(x_{2}^{0}/b_{2}).
\end{array}
\]
Hence, taking $x_{2}$ close enough to $x_{2}^{0}>0$ we get $\mathcal{L}%
(\widetilde{V})(x_{1}^{0},x_{2})>0$. So, $\widetilde{V}$ is not a viscosity
solution of (\ref{HJB}) at $(x_{1}^{0},x_{2})$. \vspace{3mm}$\Box$

\begin{remark}
In the case that $C_{\mathcal{M}}$ is empty, $A_{\mathcal{M}}=\{(0,0)\}$ and
$B_{\mathcal{M}}=M-\{(0,0)\},$ so $W(x_{2})=(\tfrac{b_{1}}{b_{2}}%
+1)x_{2}+\frac{c_{1}+c_{2}}{\lambda+q}$ and $\widetilde{V}(x_{1},x_{2}%
)=x_{1}+x_{2}+\frac{c_{1}+c_{2}}{\lambda+q}.$

This strategy is called "take the money and run". Depending on the parameters,
this strategy could be or not the optimal one.
\end{remark}

The optimal value function $V$ and the optimal strategy for surpluses in
$D^{2}\cup\mathcal{M}$, unlike the case $c_{1}/b_{1}=c_{2}/b_{2}$, cannot be
obtained in terms of the one-dimensional auxiliary optimization problem as
before. We do not have any theoretical result on the existence and structure
of the optimal strategy in $D^{2}\cup\mathcal{M}$. In Section
\ref{NumericalScheme} we will use a two-dimensional numerical scheme to
approximate the optimal strategy instead.

\section{Numerical scheme in $D^{2}\cup\mathcal{M}$\label{NumericalScheme}}

In this section we present a convergent numerical scheme to approximate the
optimal value function for the case $c_{1}/b_{1}>c_{2}/b_{2}$ in
$\mathbf{R}_{+}^{2}$. In fact, by Corollary \ref{corolario lineal en D1}, we
have that the value function satisfy%
\[
V(x_{1},x_{2})=x_{1}-\left(  b_{1}/b_{2}\right)  x_{2}+V(\left(  b_{1}%
/b_{2}\right)  x_{2},x_{2})
\]
and so it is enough to approximate the optimal value function in $D^{2}%
\cup\mathcal{M}$. This numerical scheme can be seen as a particular case of
the one described in \cite{AzMuNum} for the case that the joint multivariate
distribution function of the downward jumps is given by $F(x,y):=G(\left(
x/b_{1}\right)  \wedge\left(  y/b_{2}\right)  )$ and both the penalty and the
switch-value functions are identically zero.

Given any $\delta>0$, consider $\Delta x_{1}=c_{1}\delta$ and $\Delta
x_{2}=c_{2}\delta$, define the grid domain $\mathcal{G}^{\delta}$ in
$\mathbf{R}_{+}^{2}$ as
\[
\mathcal{G}^{\delta}:=\left\{  (n\Delta x_{1},m\Delta x_{2}):n,m\geq0\right\}
.
\]
We look, at each point of the grid $\mathcal{G}^{\delta}$, for the best local
strategy among the ones suggested by the operators of the HJB equation
(\ref{HJB}); these possible local strategies are: the first branch pays a lump
sum as dividends, the second one pays a lump sum as dividends, or none of the
branches pay dividends. These local strategies are modified in such a way that
the controlled surplus after applying these local strategies lies in the grid.
The possible control actions at any point of the grid $\mathcal{G}^{\delta}$
are defined as follows: let $(n\Delta x_{1},m\Delta x_{2})\in\mathcal{G}%
^{\delta}$ be the initial surplus and $\tau$ and $U$ be the time and size of
the first claim respectively.

\begin{enumerate}
\item $\mathbf{E}_{1}$: The first branch pays immediately $\Delta x_{1}%
=c_{1}\delta$ as dividends, so the controlled surplus becomes $((n-1)\Delta
x_{1},m\Delta x_{2})\in\mathcal{G}^{\delta}$. The control action
$\mathbf{E}_{1}$\textbf{\ }can only be applied for $n>0$.

\item $\mathbf{E}_{2}$: The second branch pays immediately $\Delta x_{2}%
=c_{2}\delta$ as dividends, so the controlled surplus becomes $(n\Delta
x_{1},(m-1)\Delta x_{2})\in\mathcal{G}^{\delta}$. The control action
$\mathbf{E}_{2}$\textbf{\ }can only be applied for $m>0$.

\item $\mathbf{E}_{0}$: Pay no dividends up to the time $\delta\wedge\tau$. In
the case that $\delta<\tau$, the uncontrolled surplus at time $\delta$ is
$((n+1)\Delta x_{1},(m+1)\Delta x_{2})\in\mathcal{G}^{\delta}$; and if
$\delta\geq\tau$, the uncontrolled surplus at time $\tau$ is $(n\Delta
x_{1}+c_{1}\tau-b_{1}U,m\Delta x_{2}+c_{2}\tau-b_{2}U)$. If this vector is in
the first quadrant, the branches pay immediately the minimum amount of
dividends in such a way that the controlled surplus lies in a point of the
grid; this end surplus can be written as
\[
(\left[  \tfrac{n\Delta x_{1}+c_{1}\tau-b_{1}U}{\Delta x_{1}}\right]  \Delta
x_{1},\left[  \tfrac{m\Delta x_{2}+c_{2}\tau-b_{2}U}{\Delta x_{2}}\right]
\Delta x_{2})\in\mathcal{G}^{\delta};
\]
the amount paid as dividends for the first and second branch at time $\tau$ is
equal to
\[
\Delta\overline{L}(\tau)=(n\Delta x_{1}+c_{1}\tau-b_{1}U-\left[
\tfrac{n\Delta x_{1}+c_{1}\tau-b_{1}U}{\Delta x_{1}}\right]  \Delta
x_{1},m\Delta x_{2}+c_{2}\tau-b_{2}U)-\left[  \tfrac{m\Delta x_{2}+c_{2}%
\tau-b_{2}U}{\Delta x_{2}}\right]  \Delta x_{2}).
\]
In the case that the end surplus $(n\Delta x_{1}+c_{1}\tau-b_{1}U,m\Delta
x_{2}+c_{2}\tau-b_{2}U)$ is not in the first quadrant, $\tau$ is the ruin time.
\end{enumerate}

For technical reasons, an extra control action $\mathbf{E}_{s}$ is considered,
under this control action no more dividends are paid. The set of control
actions is denoted by $\mathcal{E}=\{\mathbf{E}_{1},\mathbf{E}_{2}%
,\mathbf{E}_{0},\mathbf{E}_{s}\}$.

Define $\Pi_{n\Delta x_{1},m\Delta x_{2}}^{\delta}\subset\Pi_{n\Delta
x_{1},m\Delta x_{2}}$ as the set of all the admissible dividend strategies
with initial surplus $(n\Delta x_{1},m\Delta x_{2})\in\mathcal{G}^{\delta}$
which can be obtained by a sequence (finite or infinite) of control actions in
$\mathcal{E}$ at each point of the grid. The $\delta$\textit{-optimal
function} $V^{\delta}$ is defined at each point $(n\Delta x_{1},m\Delta
x_{2})\in\mathcal{G}^{\delta}$ as the supremum of the value functions of
admissible strategies which are combination of the control actions in
$\mathcal{E}$, that is
\begin{equation}
V^{\delta}(n\Delta x_{1},m\Delta x_{2})=\sup\nolimits_{\pi\in\Pi_{n\Delta
x_{1},m\Delta x_{2}}^{\delta}}V_{\pi}(n\Delta x_{1},m\Delta x_{2})\text{.}
\label{vdelta}%
\end{equation}

Note that the extra control action $\mathbf{E}_{s}$ corresponds in the
numerical scheme described in \cite{AzMuNum} to switch immediately to a
switch-value function identically $0$ (which is always smaller than $V$), so
it is never optimal to use this control; however -as we will see later on- it
is used to find a simple admissible dividends strategy in $\Pi_{n\Delta
x_{1},m\Delta x_{2}}^{\delta}$ to start a recursive construction in
$\Pi_{n\Delta x_{1},m\Delta x_{2}}^{\delta}$ which will converge to
$V^{\delta}$.

Define $v^{\delta}(n,m):=V^{\delta}(n\Delta x_{1},m\Delta x_{2})$. It is
proved in \cite{AzMuNum}, that the function $v^{\delta}$ is a solution of the
following discrete version of the HJB equation
\begin{equation}
\max\{T_{0}(W)-W,T_{1}(W)-W,T_{2}(W)-W\}=0 \label{Delta HJB}%
\end{equation}
at $(n,m)\in\mathbf{N}_{0}^{2}$. Here the operators $T_{0},$ $T_{1}$ and
$T_{2}$ are defined as%

\begin{equation}%
\begin{array}
[c]{c}%
T_{1}(W)(n,m):=W(n-1,m)+\Delta x_{1},\\
T_{2}(W)(n,m):=W(n,m-1)+\Delta x_{2},
\end{array}
\label{Definicion T1 T2}%
\end{equation}
and%

\begin{equation}
T_{0}(W)(n,m):=W(n+1,m+1)e^{-(q+\lambda)\delta}+\mathcal{I}^{\delta}(W)(n,m);
\label{Definicion T0}%
\end{equation}
where%

\begin{equation}%
\begin{array}
[c]{l}%
\mathcal{I}^{\delta}(W)(n,m)\\%
\begin{array}
[c]{ll}%
= & \int\limits_{0}^{\delta}(%
%TCIMACRO{\tint \limits_{\overline{0}}^{\frac{n\Delta x_{1}+c_{1}t}{b_{1}%
%}\wedge\frac{m\Delta x_{2}+c_{2}t}{b_{2}}}}%
%BeginExpansion
{\textstyle\int\limits_{\overline{0}}^{\frac{n\Delta x_{1}+c_{1}t}{b_{1}%
}\wedge\frac{m\Delta x_{2}+c_{2}t}{b_{2}}}}
%EndExpansion
e^{-qt}W(\left[  \frac{n\Delta x_{1}+c_{1}t-b_{1}\alpha}{\Delta x_{1}}\right]
,\left[  \frac{m\Delta x_{2}+c_{2}t-b_{2}\alpha}{\Delta x_{2}}\right]
)dG(\alpha))\lambda e^{-\lambda t}dt\\
& +\int\limits_{0}^{\delta}(%
%TCIMACRO{\tint \limits_{\overline{0}}^{\frac{n\Delta x_{1}+c_{1}t}{b_{1}%
%}\wedge\frac{m\Delta x_{2}+c_{2}t}{b_{2}}}}%
%BeginExpansion
{\textstyle\int\limits_{\overline{0}}^{\frac{n\Delta x_{1}+c_{1}t}{b_{1}%
}\wedge\frac{m\Delta x_{2}+c_{2}t}{b_{2}}}}
%EndExpansion
e^{-qt}\left(  \left(  c_{1}+c_{2}\right)  t-\alpha+n\Delta x_{1}-\left[
\frac{n\Delta x_{1}+c_{1}t-b_{1}\alpha}{\Delta x_{1}}\right]  \Delta
x_{1}\right)  dG(\alpha))\lambda e^{-\lambda t}dt\\
& +\int\limits_{0}^{\delta}(%
%TCIMACRO{\tint \limits_{\overline{0}}^{\frac{n\Delta x_{1}+c_{1}t}{b_{1}%
%}\wedge\frac{m\Delta x_{2}+c_{2}t}{b_{2}}}}%
%BeginExpansion
{\textstyle\int\limits_{\overline{0}}^{\frac{n\Delta x_{1}+c_{1}t}{b_{1}%
}\wedge\frac{m\Delta x_{2}+c_{2}t}{b_{2}}}}
%EndExpansion
e^{-qt}\left(  m\Delta x_{2}-\left[  \frac{m\Delta x_{2}+c_{2}t-b_{2}\alpha
}{\Delta x_{2}}\right]  \Delta x_{2}\right)  dG(\alpha))\lambda e^{-\lambda
t}dt.
\end{array}
\end{array}
\label{Definicion Idelta(W)}%
\end{equation}

\begin{remark}
\label{Infinitas soluciones discretas}Analogously to the Remark
\ref{Muchas soluciones viscosas} for the HJB equation (\ref{HJB}), there are
infinitely many solutions of the discrete HJB equation; in fact all the
functions $u(n,m)=n\Delta x_{1}+m\Delta x_{2}+K$ are solutions of
(\ref{Delta HJB}) for $K$ large enough. Indeed,%
\[%
\begin{array}
[c]{lll}%
(T_{0}(W)-W)(n,m) & \leq & \left(  c_{1}+c_{2}\right)  (\delta e^{-(q+\lambda
)\delta}+\frac{\lambda}{\left(  \lambda+q\right)  ^{2}}\left(  1-e^{-\left(
\lambda+q\right)  \delta}\right)  -K(1-e^{-\left(  \lambda+q\right)  \delta
})\frac{q}{\lambda+q}.
\end{array}
\]

\end{remark}

Consider the operator $T$ defined as
\begin{equation}
T:=\max\{T_{0},T_{1},T_{2}\}. \label{Definicion T}%
\end{equation}
Using the results in \cite{AzMuNum}, we have that:

\begin{itemize}
\item $v^{\delta}$ can be characterized as the smallest solution of the
discrete HJB equation (\ref{Delta HJB}). Also, if a family of strategies
$\widetilde{\pi}=\left(  \pi_{(n,m)}\right)  _{(n,m)\in\mathbf{N}_{0}^{2}}$
with $\pi_{(n,m)}\in\Pi_{n\Delta x_{1},m\Delta x_{2}}^{\delta}$ satisfy that
the function $u(n,m):=V_{\pi_{(n,m)}}(n\Delta x_{1},m\Delta x_{2})$ is a
solution of the discrete HJB equation (\ref{Delta HJB}) for all $(n,m)\in
\mathbf{N}_{0}^{2}$ then $u=v^{\delta}$ and so $V^{\delta}(n\Delta
x_{1},m\Delta x_{2})=V_{\pi_{(n,m)}}(n\Delta x_{1},m\Delta x_{2})$ for
$(n,m)\in\mathbf{N}_{0}^{2}$.

\item There exists an optimal admissible strategy for the problem
(\ref{vdelta}) at any point of the grid. This strategy, called the $\delta
$\textit{-optimal strategy}, is stationary in the sense that the control
actions depend only on which point of the grid the current surplus lies. In
the case that $T_{1}(v^{\delta})(n,m)-v^{\delta}(n,m)=0$, the optimal control
action at the point $(n\Delta x_{1},m\Delta x_{2})\in\mathcal{G}^{\delta}$ is
$\mathbf{E}_{1}$; in the case that $T_{2}(v^{\delta})(n,m)-v^{\delta}(n,m)=0,$
the optimal control action is $\mathbf{E}_{2}$; and finally in the case that
$T_{0}(v^{\delta})(n,m)-v^{\delta}(n,m)=0$, the optimal control action is
$\mathbf{E}_{0}$.

\item The function $v^{\delta}$ can be obtained recursively. The operator $T$
defined in (\ref{Definicion T}) is increasing and $v^{\delta}$ is a fixed
point of this operator. However, by Remark
\ref{Infinitas soluciones discretas}, this operator has infinitely many fixed
points. Since $T$ is not a contraction operator, $v^{\delta}$ can not be
obtained numerically as a fixed point; so we construct value functions
$v_{l}^{\delta}(n,m):=V_{\pi_{(n,m)}^{l}}(n\Delta x_{1},m\Delta x_{2})$ of
strategies $\pi_{(n,m)}^{l}\in$ $\Pi_{n\Delta x_{1},m\Delta x_{2}}^{\delta}$
which can be calculated explicitly such that $v_{l}^{\delta}$ $\nearrow$
$v^{\delta}$ as $l\rightarrow\infty$. In order to do that, let us define
iteratively families of strategies $\widetilde{\pi}_{l}=\left(  \pi
_{(n,m)}^{l}\right)  _{(n,m)\in\mathbf{N}_{0}^{2}}$ for each $l\geq1$ in the
following way;

\begin{enumerate}
\item Start with the family of strategies $\widetilde{\pi}_{1}=\left(
\pi_{(n,m)}^{1}\right)  _{(n,m)\in\mathbf{N}_{0}^{2}}$ where $\pi_{(n,m)}%
^{1}\in\Pi_{n\Delta x_{1},m\Delta x_{2}}^{\delta}$ corresponds to the local
control $\mathbf{E}_{s}$ (not more dividends are paid); the value function of
this strategy is $v_{1}^{\delta}(n,m):=V_{\pi_{(n,m)}^{1}}(n\Delta
x_{1},m\Delta x_{2})=0$.

\item Given $\widetilde{\pi}_{l}=\left(  \pi_{(n,m)}^{l}\right)
_{(n,m)\in\mathbf{N}_{0}^{2}}$ with $\pi_{(n,m)}^{l}\in\Pi_{n\Delta
x_{1},m\Delta x_{2}}^{\delta}$, define $\widetilde{\pi}_{l+1}=\left(
\pi_{(n,m)}^{l+1}\right)  _{(n,m)\in\mathbf{N}_{0}^{2}}$ as follows, for any
$(n,m)\in\mathbf{N}_{0}^{2}$, the best strategy $\pi_{(n,m)}^{l+1}\in
\Pi_{n\Delta x_{1},m\Delta x_{2}}^{\delta}$ is chosen among the ones which
initially follows one of the control actions in the set $\mathcal{E}$ and then
continues with the strategy in the family $\widetilde{\pi}_{l}$ at the end
point of the best control action. The value function of $\pi_{(n,m)}^{l+1}$ is
given by%
\[
v_{l+1}^{\delta}(n,m):=V_{\pi_{(n,m)}^{l+1}}(n\Delta x_{1},m\Delta
x_{2})=T(v_{l}^{\delta})(n,m)=T^{(l)}(v_{1}^{\delta})(n,m)~\text{for }%
(n,m)\in\mathbf{N}_{0}^{2}.
\]
Finally, since $T$ is increasing, we have that $v_{l+1}^{\delta}\geq
v_{l}^{\delta}$ for any $l\geq1$ and so there exists $\overline{v}%
=\lim_{l\rightarrow\infty}v_{l}^{\delta}$. The function $\overline{v}$ is a
solution of the discrete HJB equation (\ref{Delta HJB}) and it is constructed
as a value function of a combination of local controls in $\mathcal{E}$, then
$\overline{v}=v^{\delta}$.
\end{enumerate}

\item Let us extend the definition of $V^{\delta}$ from $\mathcal{G}^{\delta}$
to all the points in the first quadrant as
\[
V^{\delta}(x_{1},x_{2})=V^{\delta}(\left[  \frac{x_{1}}{\Delta x_{1}}\right]
\Delta x_{1},\left[  \frac{x_{2}}{\Delta x_{2}}\right]  \Delta x_{2}%
)+x_{1}-\left[  \frac{x_{1}}{\Delta x_{1}}\right]  \Delta x_{1}+x_{2}-\left[
\frac{x_{2}}{\Delta x_{2}}\right]  \Delta x_{2}.
\]

Then, $V^{\delta}(x_{1},x_{2})$ is the value function of the admissible
strategy in $\Pi_{x_{1},x_{2}}$ in which the first and second branch pay
immediately $x_{1}-\left[  x_{1}/\Delta x_{1}\right]  \Delta x_{1}$ and
$x_{2}-\left[  x_{2}/\Delta x_{2}\right]  \Delta x_{2}$ respectively as
dividends and then follows $\pi_{(n,m)}\in\Pi_{\left[  x_{1}/\Delta
x_{1}\right]  \Delta x_{1},\left[  x_{2}/\Delta x_{2}\right]  \Delta x_{2}%
}^{\delta}.$ For any $\delta>0$, it holds that $V^{\delta/2^{k}}$ $\nearrow V$
locally uniformly in the first quadrant as $k$ goes to infinity. The grids
$\mathcal{G}^{^{\delta/2^{k}}}$ are taken in order to have $\mathcal{G}%
^{^{\delta/2^{k}}}\subset$ $\mathcal{G}^{^{\delta/2^{k+1}}}$ and so
$V^{\delta/2^{k}}\leq V^{\delta/2^{k+1}}.$
\end{itemize}

\section{Numerical Examples\label{sec:examples}}

We show three numerical examples with parameters $b_{1}=b_{2}=0.5$, $c_{1}=2$,
$c_{2}=1,q=0.05$ and $\lambda=1$ and different claim-size distributions $G$.
We use the numerical scheme introduced in Section \ref{NumericalScheme} and
obtain the numerical approximation $V^{\delta}$ of the optimal value function
$V$ and the corresponding $\delta$\textit{-}optimal strategy .We can see in
general that, as we have proved in Subsection \ref{optimalStrategyIneEquality}%
, the optimal strategy for initial surplus in $D^{1}$ consist in the first
branch paying dividends immediately so the surplus moves to $D^{2}%
\cup\mathcal{M}$. It was proved in Proposition \ref{Vrulo_NeverOptimal} that
it is not optimal for some initial surplus in $\mathcal{M}$ (unlike the case
$c_{2}/c_{1}=b_{2}/b_{1}$) to pay dividends so the controlled process remains
in the set $\mathcal{M}$ until ruin time; we will see in Figures 2, 6 and
8that for some initial surpluses in $\mathcal{M}$ the best local strategy is
to get out of $\mathcal{M}$ with the first branch paying dividends immediately.

\subsection{Example 1.}

We show here a numerical example assuming an exponential claim-size
distribution $G(x)=1-e^{-dx}$ with $d=0.6$. In Figure 2 we show the $\delta
$\textit{-}optimal strategy with $\delta=0.03$, note that there is a
non-action region $\mathcal{C}^{\ast}\subset$ $D^{2}$. This figure suggests
that, as $\delta\rightarrow0$, the optimal local control in the boundary
$\mathcal{A}_{1}^{\ast}$ should be that the first branch pay some part of the
incoming premium as dividends while the second branch pay no dividends, so the
bivariate control surplus stays in $\mathcal{A}_{1}^{\ast}$ and moves
rightward to the point $(5.4,6.36)\in\mathcal{A}_{0}^{\ast}$. By contrast, the
optimal strategy in the boundary $\mathcal{A}_{2}^{\ast}$ is that none of the
two branches pay dividends so the bivariate control surplus moves inside
$\mathcal{C}^{\ast}$ because the slope of the boundary $\mathcal{A}_{2}^{\ast
}$ is larger than $c_{2}/c_{1}=1/2$.%
\[%
\begin{array}
[c]{ccc}%
%TCIMACRO{\FRAME{itbpFU}{2.2684in}{2.2987in}{0in}{\Qcb{Figure 2.}}%
%{}{estrategiaexponencial.eps}{\special{ language "Scientific Word";
%type "GRAPHIC";  maintain-aspect-ratio TRUE;  display "USEDEF";
%valid_file "F";  width 2.2684in;  height 2.2987in;  depth 0in;
%original-width 4.3284in;  original-height 4.3863in;  cropleft "0";
%filename 'EstrategiaExponencial.eps';file-properties "XNPEU";}}}%
%BeginExpansion
{\parbox[b]{2.2684in}{\begin{center}
\includegraphics[
height=2.2987in,
width=2.2684in
]%
{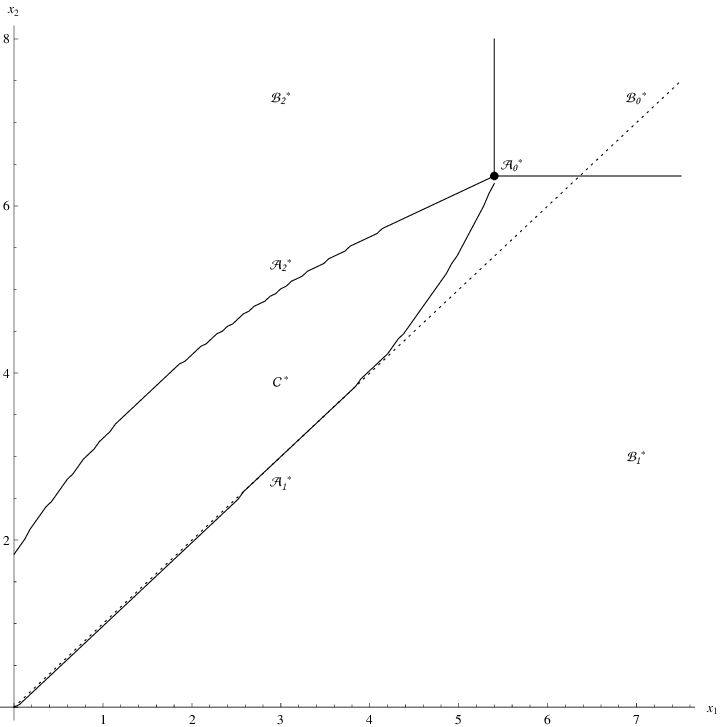}%
\\
Figure 2.
\end{center}}}%
%EndExpansion
&  &
%TCIMACRO{\FRAME{itbpFU}{2.7518in}{2.1811in}{0in}{\Qcb{Figure 3.}}%
%{}{funcionexponencial.eps}{\special{ language "Scientific Word";
%type "GRAPHIC";  maintain-aspect-ratio TRUE;  display "USEDEF";
%valid_file "F";  width 2.7518in;  height 2.1811in;  depth 0in;
%original-width 7.2195in;  original-height 5.7199in;  cropleft "0";
%filename 'FuncionExponencial.eps';file-properties "XNPEU";}}}%
%BeginExpansion
{\parbox[b]{2.7518in}{\begin{center}
\includegraphics[
%natheight=5.719900in,
%natwidth=7.219500in,
height=2.1811in,
width=2.7518in
]%
{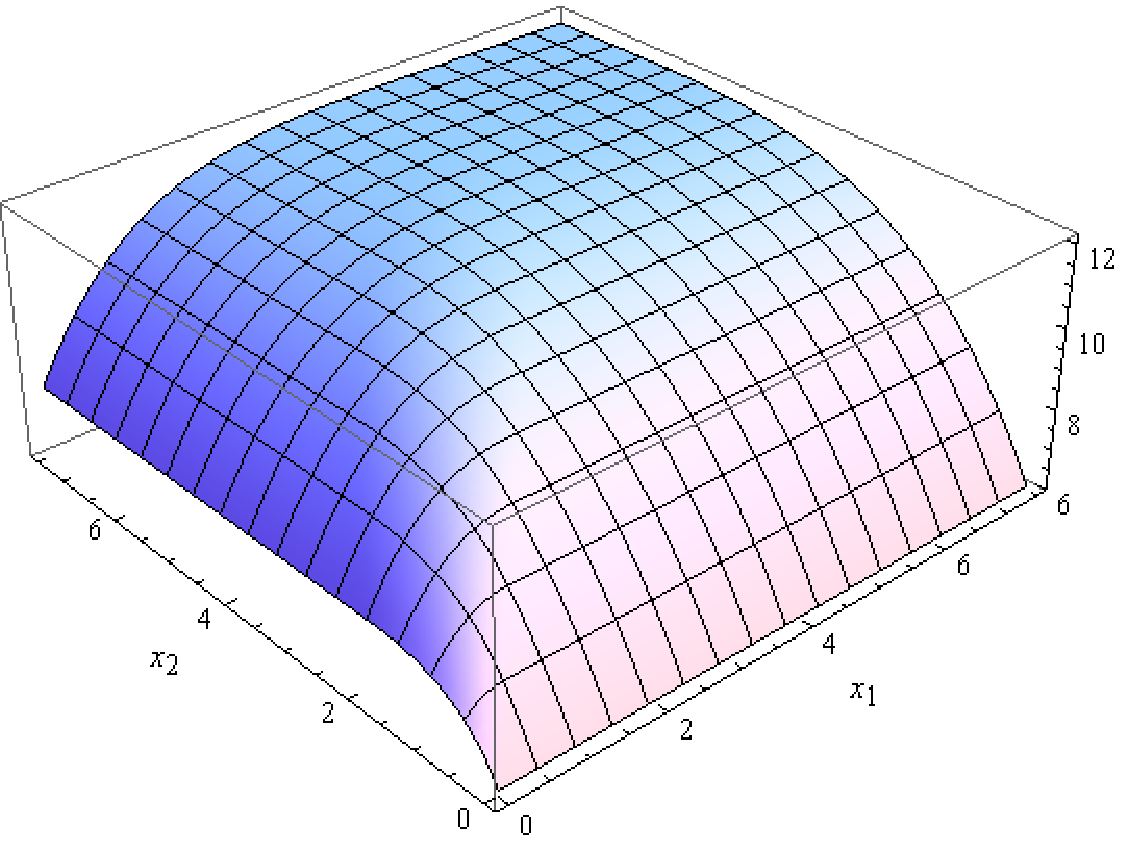}%
\\
Figure 3.
\end{center}}}%
%EndExpansion
\end{array}
\]

In Figure 3 we show $V^{\delta}$ reduced by $(x_{1}+x_{2})$. As it was noted
in Remark \ref{RemarkComparacionMerger}, the merger optimal value function
without merger cost is greater than $V$, however this could not be the case
when merger costs are considered: we compare $V^{\delta}$ with the value of
the merger optimal value function without merger cost in Figure 4, and with
the merger optimal value function with merger cost $m=3$ in Figure 5; in all
the cases the value functions are reduced by $(x_{1}+x_{2})$. We see in Figure
5 that the value function of the merger case with cost only outperforms
$V^{\delta}$ when the difference of the surpluses of the two branches is large.%

\[%
\begin{array}
[c]{ccc}%
%TCIMACRO{\FRAME{itbpFU}{3.032in}{1.887in}{0in}{\Qcb{Figure 4.}}{}%
%{optima-mergersincosto.eps}{\special{ language "Scientific Word";
%type "GRAPHIC";  maintain-aspect-ratio TRUE;  display "USEDEF";
%valid_file "F";  width 3.032in;  height 1.887in;  depth 0in;
%original-width 8.5417in;  original-height 5.3056in;  cropleft "0";
%filename 'Optima-MergerSinCosto.eps';file-properties "XNPEU";}}}%
%BeginExpansion
{\parbox[b]{3.032in}{\begin{center}
\includegraphics[
%natheight=5.305600in,
%natwidth=8.541700in,
height=1.887in,
width=3.032in
]%
{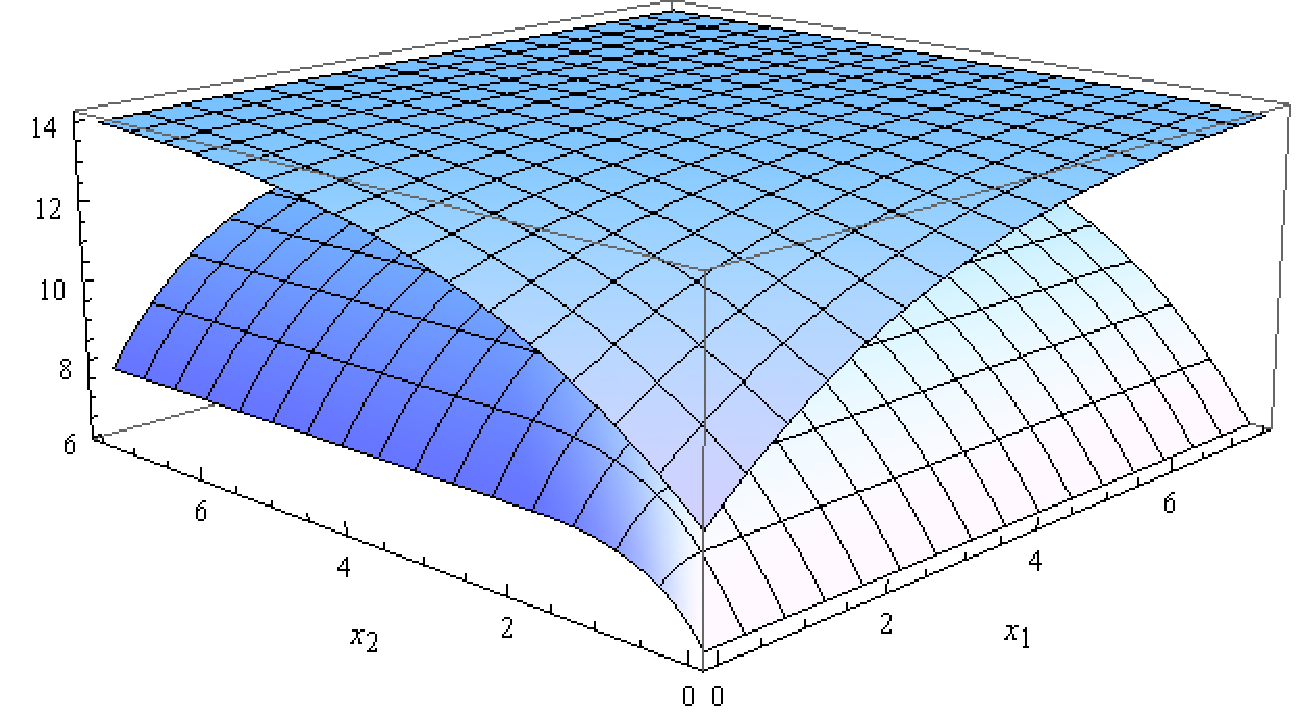}%
\\
Figure 4.
\end{center}}}%
%EndExpansion
&  &
%TCIMACRO{\FRAME{itbpFU}{2.7501in}{1.7953in}{0in}{\Qcb{Figure 5.}}%
%{}{optima-mergercosto.eps}{\special{ language "Scientific Word";
%type "GRAPHIC";  maintain-aspect-ratio TRUE;  display "USEDEF";
%valid_file "F";  width 2.7501in;  height 1.7953in;  depth 0in;
%original-width 7.6942in;  original-height 5.0142in;  cropleft "0";
%filename 'Optima-MergerCosto.eps';file-properties "XNPEU";}}}%
%BeginExpansion
{\parbox[b]{2.7501in}{\begin{center}
\includegraphics[
%natheight=5.014200in,
%natwidth=7.694200in,
height=1.7953in,
width=2.7501in
]%
{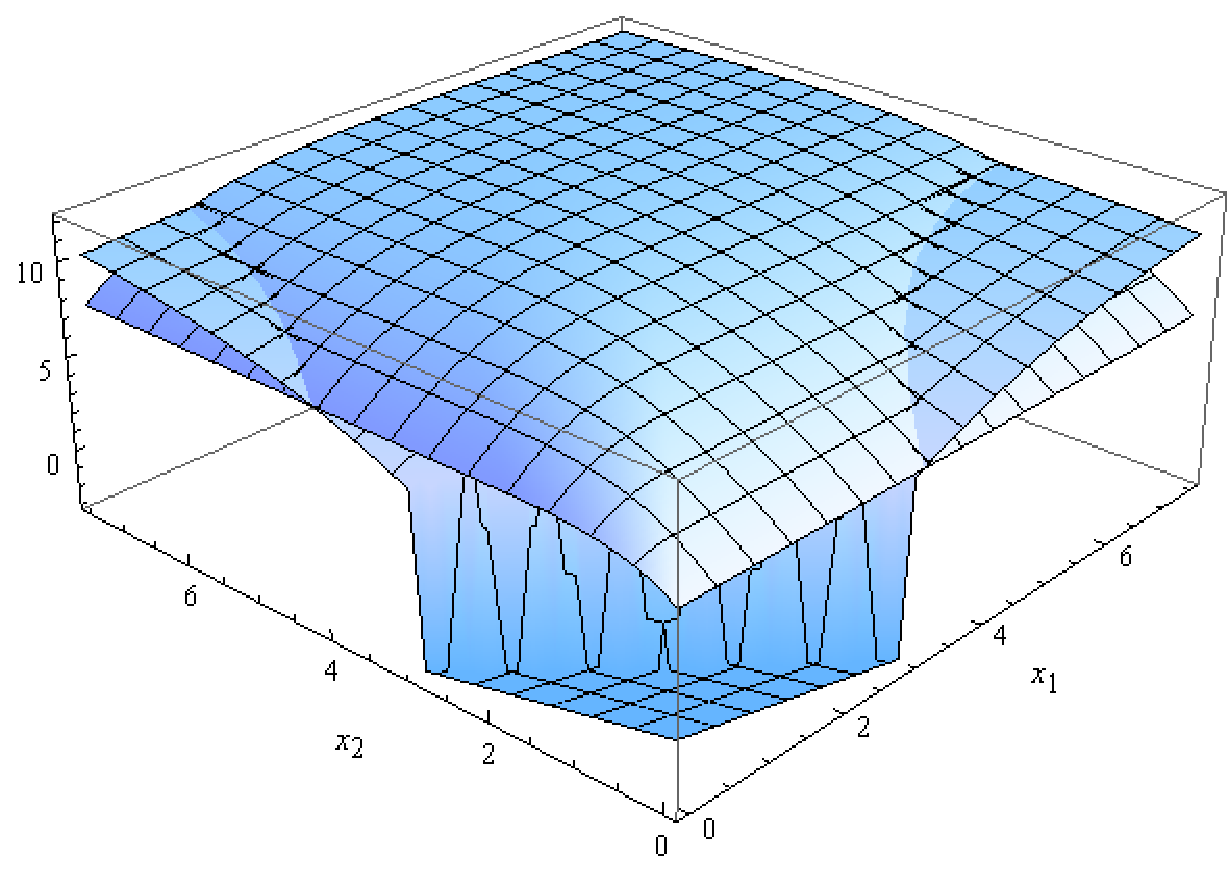}%
\\
Figure 5.
\end{center}}}%
%EndExpansion
\end{array}
\]

\subsection{Example 2.}

We consider here a gamma claim-size distribution
\[
G(x)=1-(1+\frac{6}{7}x)e^{-\frac{6}{7}x}.
\]
We show in Figures 6 and 7 the $\delta$\textit{-}optimal strategy and
$V^{\delta}$ reduced by $(x_{1}+x_{2})$ respectively for $\delta=0.025$. Note
that, unlike the previous example, the sets $\mathcal{A}_{0}^{\ast
}=\{(0,0),(4.00,4.75)\}$ and $\mathcal{B}_{0}^{\ast}$ have two connected
components. The graph suggest that $V^{\delta}$ is not differentiable at the
lower boundary between $\mathcal{B}_{1}^{\ast}$ and $\mathcal{B}_{0}^{\ast}$
but it is differentiable at the upper one (this mirrors the smooth fitting
property in the one-dimensional case); a similar observation could be made
about the boundaries between $\mathcal{B}_{2}^{\ast}$ and $\mathcal{B}%
_{0}^{\ast}.$%

\[%
\begin{array}
[c]{ccc}%
%TCIMACRO{\FRAME{itbpFU}{2.2917in}{2.3203in}{0in}{\Qcb{Figure 6.}}%
%{}{estrategiagamma.eps}{\special{ language "Scientific Word";
%type "GRAPHIC";  maintain-aspect-ratio TRUE;  display "USEDEF";
%valid_file "F";  width 2.2917in;  height 2.3203in;  depth 0in;
%original-width 4.8421in;  original-height 4.9018in;  cropleft "0";
%filename 'EstrategiaGamma.eps';file-properties "XNPEU";}}}%
%BeginExpansion
{\parbox[b]{2.2917in}{\begin{center}
\includegraphics[
height=2.3203in,
width=2.2917in
]%
{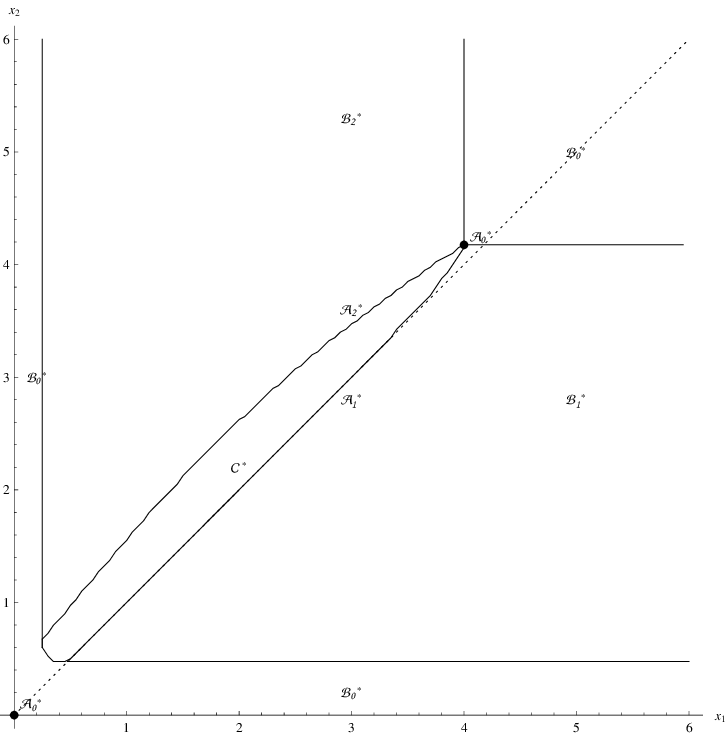}%
\\
Figure 6.
\end{center}}}%
%EndExpansion
&  &
%TCIMACRO{\FRAME{itbpFU}{3.1548in}{2.0574in}{0in}{\Qcb{Figure 7.}}%
%{}{funciongamma.eps}{\special{ language "Scientific Word";  type "GRAPHIC";
%width 3.1548in;  height 2.0574in;  depth 0in;  original-width 7.958in;
%original-height 5.1802in;  cropleft "0";  croptop "1";  cropright "1";
%cropbottom "0";  filename 'FuncionGamma.eps';file-properties "XNPEU";}}}%
%BeginExpansion
{\parbox[b]{3.1548in}{\begin{center}
\includegraphics[
%natheight=5.180200in,
%natwidth=7.958000in,
height=2.0574in,
width=3.1548in
]%
{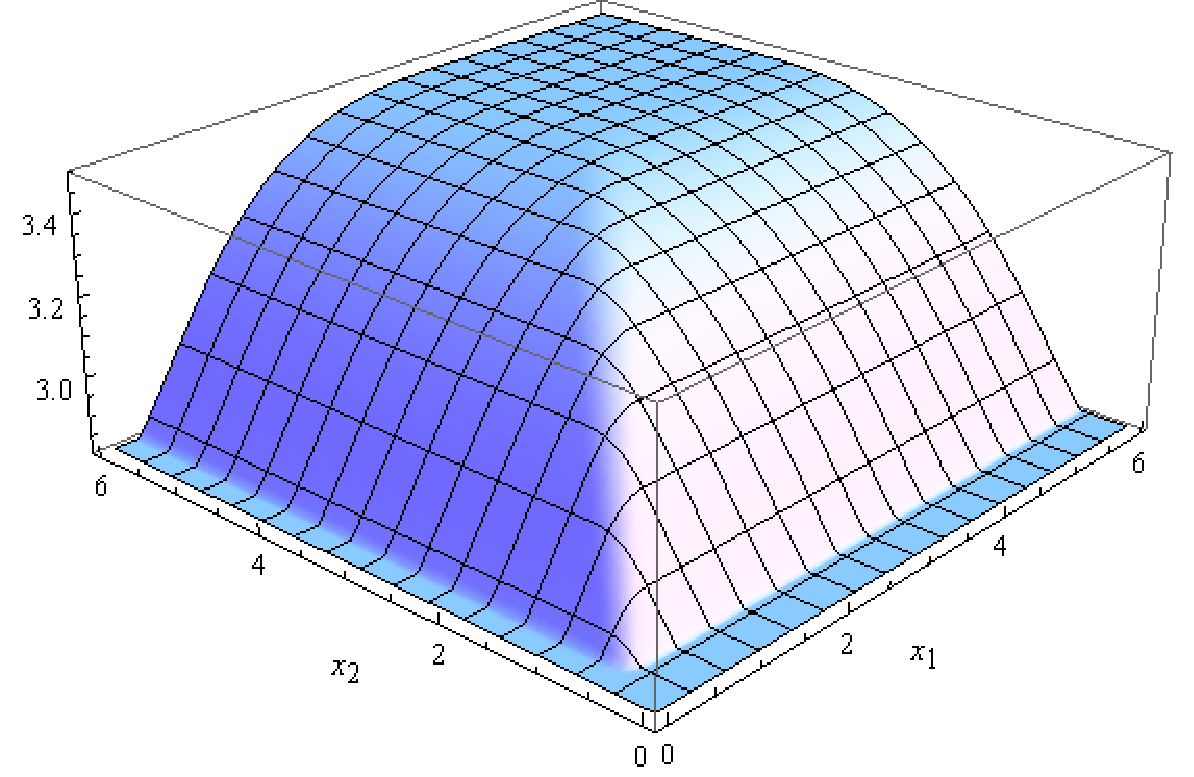}%
\\
Figure 7.
\end{center}}}%
%EndExpansion
\end{array}
\]

As in the first example, Figure 6 suggests that the optimal local control in
the boundary $\mathcal{A}_{1}^{\ast}$ should be that the first branch pay some
part of the incoming premium as dividends while the second branch pay no
dividends, so the bivariate control surplus stays in $\mathcal{A}_{1}^{\ast}$
and moves rightward to the point $(4.00,4.75)\in\mathcal{A}_{0}^{\ast}$.
Again, the optimal strategy in the boundary $\mathcal{A}_{2}^{\ast}$ is that
none of the two branches pay dividends so the bivariate control surplus moves
inside $\mathcal{C}^{\ast}$ because the slope of the boundary $\mathcal{A}%
_{2}^{\ast}$ is larger than $c_{2}/c_{1}=1/2$.

\subsection{Example 3.}

Finally, we consider claims with constant size $\alpha_{0}=29/12.$ Figures 8
and 9 show the $\delta$\textit{-}optimal strategy and $V^{\delta}$ reduced by
$(x_{1}+x_{2})$ respectively for $\delta=0.02$. As in the previous example,
the sets $\mathcal{A}_{0}^{\ast}=\{(0,0),(3.56,3.62)\}$ and $\mathcal{B}%
_{0}^{\ast}$ have two connected components.%

\[%
\begin{array}
[c]{ccc}%
%TCIMACRO{\FRAME{itbpFU}{2.4569in}{2.4889in}{0in}{\Qcb{Figure 8.}}%
%{}{estrategiasalto.eps}{\special{ language "Scientific Word";
%type "GRAPHIC";  maintain-aspect-ratio TRUE;  display "USEDEF";
%valid_file "F";  width 2.4569in;  height 2.4889in;  depth 0in;
%original-width 4.5645in;  original-height 4.6241in;  cropleft "0";
%filename 'EstrategiaSalto.eps';file-properties "XNPEU";}}}%
%BeginExpansion
{\parbox[b]{2.4569in}{\begin{center}
\includegraphics[
height=2.4889in,
width=2.4569in
]%
{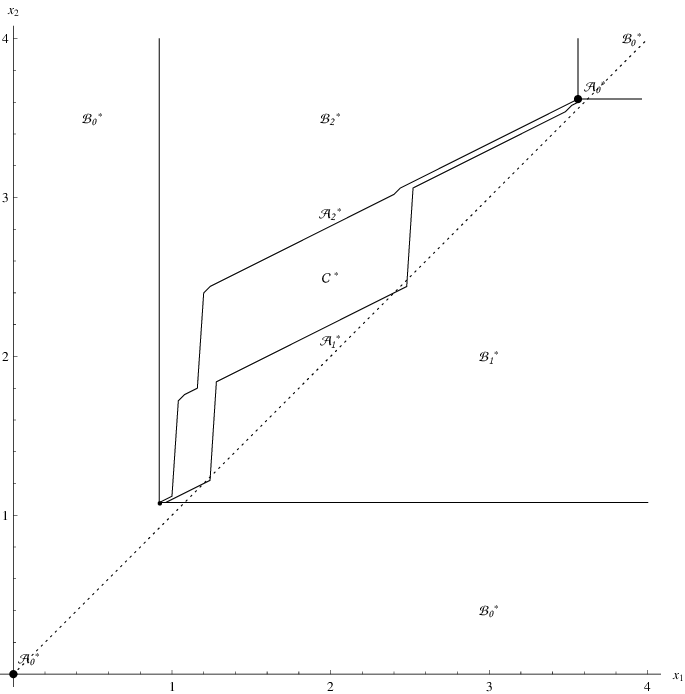}%
\\
Figure 8.
\end{center}}}%
%EndExpansion
&  &
%TCIMACRO{\FRAME{itbpFU}{3.2595in}{2.4076in}{0in}{\Qcb{Figure 9.}}%
%{}{funcionsalto.eps}{\special{ language "Scientific Word";  type "GRAPHIC";
%width 3.2595in;  height 2.4076in;  depth 0in;  original-width 8.1526in;
%original-height 6.0139in;  cropleft "0";  croptop "1";  cropright "1";
%cropbottom "0";  filename 'FuncionSalto.eps';file-properties "XNPEU";}}}%
%BeginExpansion
{\parbox[b]{3.2595in}{\begin{center}
\includegraphics[
%natheight=6.013900in,
%natwidth=8.152600in,
height=2.4076in,
width=3.2595in
]%
{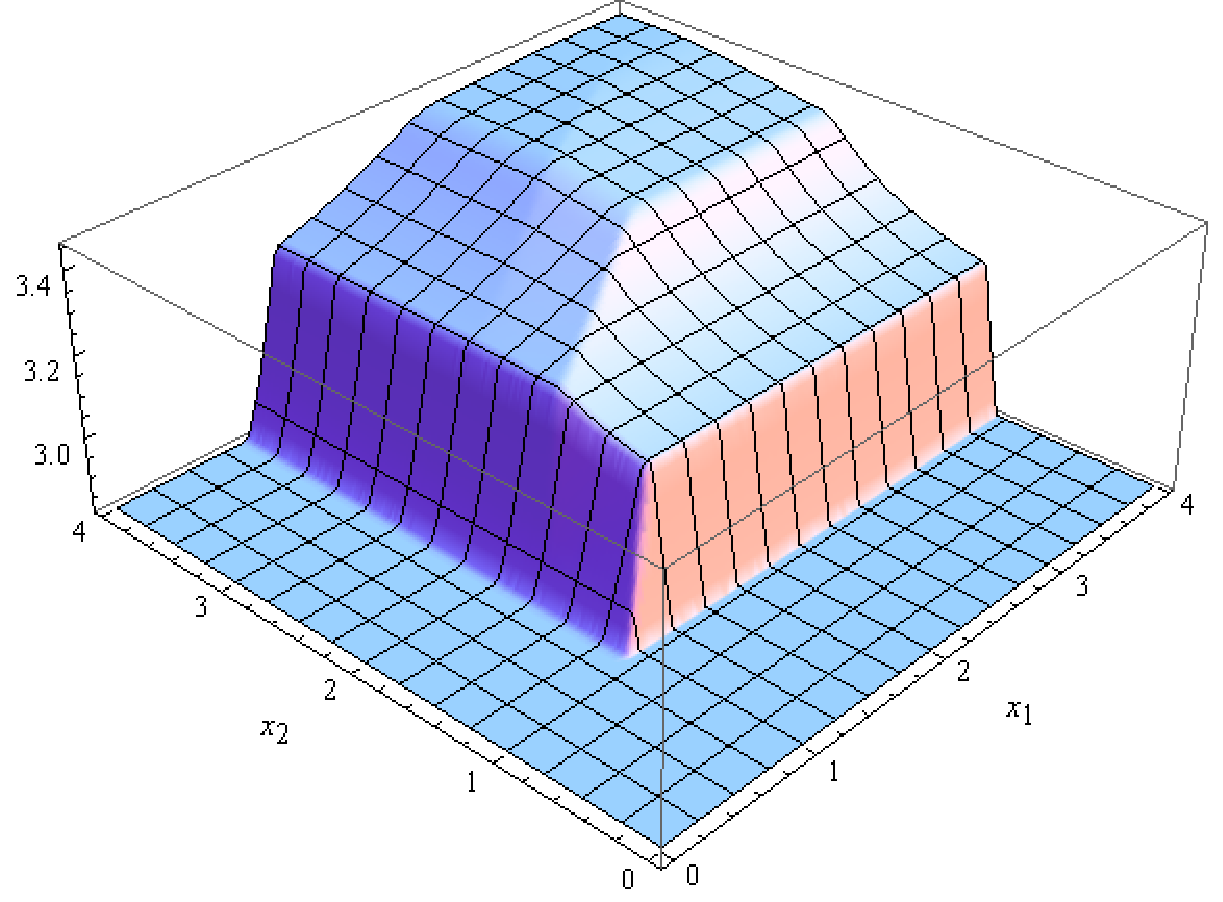}%
\\
Figure 9.
\end{center}}}%
%EndExpansion
\end{array}
\]

Note that there is a relation between the shape of $\mathcal{A}_{1}^{\ast}$
and $\mathcal{A}_{2}^{\ast}$ , the constant size of the claims $U=29/12$ and
the rate of growth $c_{2}/c_{1}=1/2$ of the uncontrolled bivariate surplus in
the event of no claims: $\mathcal{A}_{1}^{\ast}$ and $\mathcal{A}_{2}^{\ast}$
contain segments with slope $c_{2}/c_{1}$ and $\Delta x_{1}=b_{1}U$. As in the
previous examples, the optimal control strategy en $\mathcal{A}_{2}^{\ast}$ is
that none of the two branches pay dividends, and also the optimal strategy in
the boundary $\mathcal{A}_{1}^{\ast}$ consists on the second branch paying no
dividends and the first branch paying some part of the incoming premium as
dividends in such a way that the bivariate control surplus stays in
$\mathcal{A}_{1}^{\ast}$ and moves rightward to the point $(3.56,3.62)\in
\mathcal{A}_{0}^{\ast}$. In the segments of $\mathcal{A}_{1}^{\ast}$ with
slope $c_{2}/c_{1}$, the first branch does not need to pay dividends in order
to remain in $\mathcal{A}_{1}^{\ast}$.

\section*{Acknowledgement}

This research is support by the FP7 Grant PIRSES-GA-2012-318984. Zbigniew
Palmowski is supported by the Ministry of Science and Higher Education of
Poland under the grant 2013/09/B/HS4/01496.

\end{document}